\documentclass[aos,seceqn,nameyear,dvips]{arximspdf}
\usepackage{accents}
\usepackage{dcolumn}
\usepackage{graphicx}

\doi{10.1214/10-AOS824}
\volume{38}
\issue{6}
\pubyear{2010}
\firstpage{3630}
\lastpage{3659}

\makeatletter

\newcolumntype{d}[1]{D{.}{.}{#1}}

\newtheorem{theorem}{Theorem}
\newproclaim{remark}{Remark}

\makeatother

\begin{document}
\begin{frontmatter}

\title{ANOVA for longitudinal data with missing values\thanksref{T1}}
\runtitle{ANOVA for longitudinal data}

\thankstext{T1}{Supported by NSF Grants SES-0518904, DMS-06-04563 and DMS-07-14978.}

\begin{aug}
\author[A]{\fnms{Song Xi} \snm{Chen}\corref{}\ead[label=e1]{songchen@iastate.edu}} and
\author[B]{\fnms{Ping-Shou} \snm{Zhong}\ead[label=e3]{pszhong@iastate.edu}}
\runauthor{S. X. Chen and P.-S. Zhong}
\affiliation{Iowa State University and Peking University, and Iowa
State University}
\address[A]{Department of Statistics\\
Iowa State University\\
Ames, Iowa 50011-1210\\
USA\\
and\\
Guanghua School of Management\\
Center for Statistical Science\\
Peking University, Beijing 100871\\
China\\
\printead{e1}}
\address[B]{Department of Statistics\\
Iowa State University\\
Ames, Iowa 50011-1210\\
USA\\
\printead{e3}}
\end{aug}

\received{\smonth{2} \syear{2010}}
\revised{\smonth{4} \syear{2010}}

%
\begin{abstract}
We carry out ANOVA comparisons of multiple treatments for longitudinal
studies with missing values. The treatment effects are modeled
semiparametrically via a partially linear regression which is flexible
in quantifying the time effects of treatments. The empirical likelihood
is employed to formulate model-robust nonparametric ANOVA tests for
treatment effects with respect to covariates, the nonparametric
time-effect functions and interactions between covariates and time. The
proposed tests can be readily modified for a variety of data and model
combinations, that encompasses parametric, semiparametric and
nonparametric regression models; cross-sectional and longitudinal data,
and with or without missing values.
\end{abstract}

%
\begin{keyword}[class=AMS]
\kwd[Primary ]{62G10}
\kwd[; secondary ]{62G20}
\kwd{62G09}.
\end{keyword}
\begin{keyword}
\kwd{Analysis of variance}
\kwd{empirical likelihood}
\kwd{kernel smoothing}
\kwd{missing at random}
\kwd{semiparametric model}
\kwd{treatment effects}.
\end{keyword}

\end{frontmatter}

\section{Introduction}\label{sec1}
Randomized clinical trials and observational studies are often
used to evaluate treatment effects. While the treatment
versus control studies are popular, multi-treatment comparisons
beyond two samples are commonly practised in clinical trails
and observational studies.
In addition to evaluate overall treatment effects, investigators
are also interested in intra-individual changes over time by
collecting repeated measurements on each individual over time.
Although most longitudinal studies are desired to have all
subjects measured at the same set of time points, such ``balanced''
data may not be available in practice due to missing values.
Missing values arise when scheduled measurements are not made,
which make the data ``unbalanced.'' There is a good body of
literature on parametric, nonparametric and semiparametric
estimation for longitudinal data with or without missing values. This
includes \citet{r18}, \citet{r16}, \citet{r43},
\citet{r44}, \citet{r9} for methods developed
for longitudinal data without missing values; and \citet{r22},
\citet{r23}, \citet{r15}, \citet{r32} for missing values.

The aim of this paper is to develop ANOVA tests for
multi-treatment comparisons in longitudinal studies with or
without missing values. Suppose that at time~$t$, corresponding to
$k$ treatments there are $k$ mutually independent samples,
%
\[
\{(Y_{1i}(t),X^\tau_{1i}(t)) \}_{i=1}^{n_1},\ldots,
\{(Y_{ki}(t),X^\tau_{ki}(t))\}_{i=1}^{n_k},
\]
where
the response variable $Y_{ji}(t)$ and the covariate $X_{ji}(t)$
are supposed to be measured at time points $t=t_{ji
1},\ldots,t_{jiT_j}$. Here $T_j$ is the fixed number of scheduled
observations for the $j$th treatment. However,
$\{Y_{ji}(t),X^\tau_{ji}(t)\}$ may not be observed at some times,
resulting in missing values in either the response $Y_{ji}(t)$
or the covariates $X_{ji}(t)$.

We consider a semiparametric regression
model for the longitudinal data
%
%
\begin{eqnarray} \label{eq:model0}
Y_{ji}(t)=
X_{ji}^{\tau}(t)\beta_{j0}+M^\tau(X_{ji}(t),t)\gamma
_{j0}+g_{j0}(t)+\varepsilon_{ji}(t),\nonumber\\[-10pt]\\[-10pt]
\eqntext{j=1,2,\ldots,k,}
\end{eqnarray}
where $M(X_{ji}(t),t)$ are known functions of $X_{ji}(t)$ and time
$t$ representing interactions between the covariates and the time,
$\beta_{j0}$ and $\gamma_{j0}$ are $p$- and $q$-dimensional
parameters,
respectively, $g_{j0}(t)$ are unknown smooth functions representing the time
effect, and $\{\varepsilon_{ji}(t)\}$ are residual time series.
Such a semiparametric model may be viewed as an extended partially linear
model.
The partially linear model has been used for
longitudinal data analysis; see \citet{r48}, \citet{r49},
\citet{r19}, \citet{r39}. \citet{r43} and \citet{r44}
proposed estimation and confidence regions for a
semiparametric varying coefficient regression model.
Despite a body of works on estimation for longitudinal data,
analysis of variance for longitudinal data have attracted much
less attention. A few exceptions include \citet{r10} who
proposed an ANOVA test in a fully parametric setting; and \citet
{r34} who considered a two sample test in a fully
nonparametric setting.

In this paper, we propose ANOVA tests for differences among the
$\beta_{j0}$'s and the baseline time functions $g_{j0}$'s,
respectively, in the presence of the interactions. The ANOVA
statistics are formulated based on the empirical likelihood [Owen
(\citeyear{r26}, \citeyear{r28})],
which can be viewed as a nonparametric
counterpart of the conventional parametric likelihood.
Despite its not requiring a fully parametric model, the empirical
likelihood
enjoys two key properties of a conventional likelihood, the Wilks'
theorem [\citet{r27}, \citet{r30}, \citet{r8}]
and Bartlett correction [\citet{r5}, \citet{r1}];
see \citet{r3a} for an overview on the
empirical likelihood for regression. This resemblance to the
parametric likelihood
ratio motivates us to consider\vadjust{\goodbreak} using empirical likelihood to formulate
ANOVA test
for longitudinal data in nonparametric situations. This will
introduce a much needed model-robustness in the ANOVA testing.

Empirical likelihood has been used in studies for either missing
or longitudinal data. \citet{r41}, \citet{r40} considered
an empirical likelihood inference with a kernel regression
imputation for missing responses. \citet{r17} treated
estimation for the partially linear model with missing covariates.
For longitudinal data, Xue and Zhu (\citeyear{r45}, \citeyear{r46})
proposed a bias
correction method to make the empirical likelihood statistic asymptotically
pivotal in a one sample partially linear model; see also You, Chen
and Zhou (\citeyear{r47}) and \citet{r14}.

In this paper, we propose three empirical likelihood based ANOVA tests
for the
equivalence of the treatment effects with respect to (i) the
covariate $X_{ji}$; (ii) the interactions $M(X_{ji}(t),t)$ and
(iii) the time effect functions $g_{j0}(\cdot)$'s, by formulating
empirical likelihood ratio test statistics. It is shown that for the proposed
ANOVA tests for the covariates effects and the interactions, the
empirical likelihood ratio statistics are asymptotically chi-squared
distributed,
which resembles the conventional ANOVA statistics based on
parametric likelihood ratios. This is achieved without parametric
model assumptions for the residuals in the presence of the
nonparametric time effect functions and missing values. Hence, the
empirical likelihood ANOVA tests have the needed model-robustness. Another
attraction of the proposed ANOVA tests is that they encompass a
set of ANOVA tests for a variety of data and model combinations.
Specifically, they imply specific ANOVA tests for both
cross-sectional and longitudinal data; for parametric,
semiparametric and nonparametric regression models; and with or
without missing values.
%


The paper is organized as below. In Section~\ref{sec2}, we describe the
model and the missing value mechanism. Section~\ref{sec3} outlines the
ANOVA test for comparing treatment effects due to the covariates:
whereas the tests regarding interaction are proposed
in Section~\ref{sec5}. Section~\ref{sec4} considers ANOVA test for the
nonparametric time
effects. The bootstrap
calibration to the ANOVA test on the nonparametric part is
outlined in Section~\ref{sec6}. Section~\ref{sec7} reports simulation
results. We
applied the proposed ANOVA tests in Section~\ref{sec8} to analyze an HIV-CD4
data set.
Technical assumptions are presented in the
\hyperref[app]{Appendix}. All the technical proofs to the theorems are
reported in a
supplement article [\citet{r4}].

\section{Models, hypotheses and missing values}\label{sec2}

For the $i$th individual of the $j$th treatment, the measurements
taken at time $t_{jim}$ follow a semiparametric model
%
%
\begin{eqnarray}\label{eq:model}
Y_{ji}(t_{jim})&=&X^{\tau}_{ji}(t_{jim})\beta_{j0}+M^\tau
(X_{ji}(t_{jim}),t_{jim})\gamma_{j0}\nonumber\\[-8pt]\\[-8pt]
&&{}+g_{j0}(t_{jim})+\varepsilon
_{ji}(t_{jim}),\nonumber
\end{eqnarray}
for $j=1,\ldots,k$, $i =1,\ldots,n_j$, $m=1,
\ldots, T_j$.
Here $\beta_{j0}$ and $\gamma_{j0}$ are unknown $p$- and
$q$-dimensional parameters\vspace*{-2pt}
and $g_{j0}(t)$ are unknown functions
representing the time effects of the treatments. The time points $\{
t_{jim}\}_{m=1}^{T_j}$ 
are known design points. For   ease of notation, we write
$(Y_{jim},X_{jim}^\tau,M_{jim}^\tau)$ to denote
$(Y_{ji}(t_{jim}),X_{ji}^\tau(t_{jim}),M^\tau(X_{ji}(t_{jim}),t_{jim}))$.
Also, we will use
$\mathbb{X}_{jim}^\tau=(X_{jim}^\tau,M_{jim}^\tau)$ and
$\mathbb{\xi}_j^\tau=(\beta_j^\tau,\gamma_j^\tau)$. For each
individual, the residuals $\{\epsilon_{ji}(t)\}$ satisfy $E\{
\varepsilon_{ji}(t)|X_{ji}(t)\}=0$, $\operatorname{Var}\{
\varepsilon_{ji}(t)|X_{ji}(t)\}=\sigma^2_{j}(t)$ and
\[
\operatorname{Cov}\{ \varepsilon_{ji}(t), \varepsilon_{ji}(s) | X_{ji}(t),
X_{ji}(s) \} = \rho_{j}(s,t) \sigma_{j}(t) \sigma_{j}(s),
\]
where $\rho_{j}(s,t)$ is the conditional correlation coefficient
between two residuals at two different times. And the residual
time series $\{\epsilon_{ji}(t)\}$ from different subjects and
different treatments are independent. Without loss of generality,
we assume $t, s \in[0,1]$. For the purpose of identifying
$\beta_{j0}$, $\gamma_{j0}$ and $g_{j0}(t)$, we assume
%
\[
(\beta_{j0},\gamma_{j0},g_{j0})=\mathop{\arg\min}_{(\beta_{j},\gamma
_{j},g_{j})}\frac{1}{n_jT_j}\sum_{i=1}^{n_j} \sum_{m=1}^{T_j}
E\{Y_{jim}-X^{\tau}_{jim}\beta_{j}-M^\tau_{jim}\gamma
_{j}-g_{j}(t_{jim})\}^2.
\]
We also require that
$\frac{1}{n_jT_j}\sum_{i=1}^{n_j}\sum_{m=1}^{T_j}E(\widetilde
{\mathbb{X}}_{jim}\widetilde{\mathbb{X}}_{jim}^\tau)>0$,
where
$\widetilde{\mathbb{X}}_{jim}=\mathbb{X}_{jim}-E(\mathbb{X}_{jim}|t_{jim})$.
This condition also rules out $M(X_{ji}(t), t)$ being a pure
function of~$t$, and hence it has to be genuine interaction. For
the same reason, the intercept in model (\ref{eq:model}) is
absorbed into the nonparametric part $g_{j0}(t)$.

As commonly exercised in the
partially linear model [\citet{r37}; \citet{r20}],
there is a
secondary model for the covariate~$X_{jim}$:
%
%
\begin{eqnarray}
\label{Secondmodel}
X_{jim}=h_j (t_{jim})+u_{jim},\nonumber\\[-8pt]\\[-8pt]
\eqntext{j=1,2,\ldots,k, i=1,\ldots,n_j, m=1,\ldots,T_j,}
\end{eqnarray}
where $h_j(\cdot)$'s are $p$-dimensional smooth functions with
continuous second derivatives, the residual\vspace*{1pt}
$u_{jim}=(u_{jim}^1,\ldots,u_{jim}^p)^\tau$ satisfy $E(u_{jim})=0$
and $u_{jl}$ and $u_{jk}$ are independent for $l \ne k$, where
$u_{jl}=(u_{jl1},\ldots,u_{jlT_j})$. By the identification
condition given above, the covariance matrix of $u_{jim}$ is
assumed to be finite and positive definite.

We are interested in testing three ANOVA hypotheses. The first one is on
the treatment effects with respect to the covariates:
\[
H_{0a}\dvtx\beta_{10}=\beta_{20}=\cdots=\beta_{k0} \quad\mbox{vs.}\quad
H_{1a}\dvtx\beta_{i 0} \ne\beta_{j 0}\qquad\mbox{for some } i \ne
j.
\]
The second one is regarding the time effect functions:
\[
H_{0b}\dvtx
g_{10}(\cdot)=\cdots=g_{k0}(\cdot) \quad\mbox{vs.}\quad H_{1
b}\dvtx g_{i 0}(\cdot) \ne g_{j 0}(\cdot)\qquad\mbox{for some }i \ne j.
\]
The third one is on the existence of the interaction
$H_{0c}\dvtx\gamma_{j0}=0$ and $H_{1c}\dvtx\break\gamma_{j0}\neq0$.
And the last one is the ANOVA test for
\[
H_{0d}\dvtx
\gamma_{10}=\gamma_{20}=\cdots=\gamma_{k0} \quad\mbox{vs.}\quad
H_{1d}\dvtx\gamma_{i0}\neq\gamma_{j0}\qquad\mbox{for some }i\neq j.
\]

Let $X_{ji}=\{X_{ji0},\ldots,X_{jiT_j}\}$ and
$Y_{ji}=\{Y_{ji0},\ldots,Y_{jiT_j}\}$ be the complete time series
of the covariates and responses of the $(j,i)$th subject (the $i$th
subject in the $j$th treatment), and
$\accentset{\leftharpoonup}{Y}_{jit,d}=\{Y_{ji (t-d)},\ldots
,Y_{ji(t-1)}\}$
and $\accentset{\leftharpoonup}{X}_{jit,d}=\{X_{ji
(t-d)},\ldots,X_{ji(t-1)}\}$ be the past $d$ observations at time
$t$ for a positive integer $d \le\min_j \{T_j\}$. For $t<d$, we set $d=t-1$.

Define the missing value indicator $\delta_{jit}=1$ if
$(X^\tau_{ji t}, Y_{ji t})$ is observed and $\delta_{jit}=0$ if
$(X^\tau_{ji t}, Y_{ji t})$ is missing. Here, we assume $X_{ji t}$
and $Y_{jit}$ are either both observed or both missing. This
simultaneous missingness of $X_{ji t}$ and $Y_{jit}$ is
for the ease of mathematical exposition.
We also assume that $\delta_{ji0}=1$, namely the first visit of
each subject is always made.

Monotone missingness is a common assumption in the analysis of
longitudinal data [\citet{r32}]. It assumes that
if $\delta_{ji (t-1)}=0$ then $\delta_{jit}=0$. However, in
practice after missing some scheduled appointments people may
rejoin the study. This kind of casual drop-out appears quite
often in empirical studies.
To allow more data being included in the analysis, we relax the
monotone missingness to allow segments of consecutive $d$ visits
being used. Let $\delta_{jit,d} = \prod_{l=1}^{d} \delta_{ji
(t-l)}$. We assume the missingness of $(X^\tau_{jit}, Y_{jit})$
is missing at random (MAR) \citet{r33} given its immediate past
$d$ complete observations, namely
%
%
\begin{eqnarray}\label{eq:MAR}\qquad
P(\delta_{jit}=1|\delta_{jit,d}=1, X_{ji},Y_{ji})
&=& P(\delta_{jit}=1|\delta_{jit,d}=1,\accentset{\leftharpoonup
}{X}_{jit,d},\accentset{\leftharpoonup}{Y}_{jit,d})
\nonumber\\[-8pt]\\[-8pt]
&=& p_{j}(\accentset{\leftharpoonup}{X}_{jit,d},
\accentset{\leftharpoonup}{Y}_{jit,d};\theta_{j0}).\nonumber
\end{eqnarray}
Here the missing propensity $p_{j}$ is known up to a
parameter $\theta_{j0}$.
To allow derivation of a binary likelihood function, we need to
set $\delta_{jit}=0$ if $\delta_{jit,d}=0$ when there is some
drop-outs among the past $d$ visits, which is only temporarily if
$\delta_{jit} =1$. This set-up ensures
%
%
\begin{equation}\label{eq:force}
P(\delta_{jit}=0|\delta_{jit,d}=0,\accentset{\leftharpoonup
}{X}_{jit,d},\accentset{\leftharpoonup}{Y}_{jit,d})=1.
\end{equation}

Now the conditional binary likelihood for
$\{\delta_{jit}\}_{t=1}^{T_j}$ given $X_{ji}$ and $Y_{ji}$ is
%
\begin{eqnarray*}
&& P(\delta_{ji0},\ldots,\delta_{jiT_j}|X_{ji},Y_{ji})\\
&&\qquad=\prod_{m=1}^{T_j}P\bigl(\delta_{jim}|\delta_{ji(m-1)},\ldots,\delta
_{ji0},X_{ji},Y_{ji}\bigr)\\
&&\qquad=\prod_{m=1}^{T_j}P(\delta_{jim}|\delta_{jim,d}=1,\accentset
{\leftharpoonup}{X}_{jim,d},\accentset{\leftharpoonup
}{Y}_{jim,d})\\
&&\qquad=\prod_{m=1}^{T_j}\bigl[
p_{j}(\accentset{\leftharpoonup}{X}_{jim,d},\accentset{\leftharpoonup
}{Y}_{jim,d};\theta_j)^{\delta_{jim}}
\{1-p_{j}(\accentset{\leftharpoonup}{X}_{jim,d},\accentset
{\leftharpoonup}{Y}_{jim,d};\theta_{j})\}^{(1-\delta_{jim})}\bigr]^{\delta_{jim,d}}.
\end{eqnarray*}
In the second equation above, we use both the MAR in (\ref{eq:MAR})
and (\ref{eq:force}).
Hence, the parameters $\theta_{j 0}$ can
be estimated by maximizing the binary likelihood
%
%
\begin{eqnarray}
\label{BinLikehood} {\mathcal{L}}_{B_j}(\theta_j)
&=&\prod_{i=1}^{n_{j}}\prod_{t=1}^{T_j}\bigl[
p_{j}(\accentset{\leftharpoonup}{X}_{jit,d},\accentset{\leftharpoonup
}{Y}_{jit,d};\theta_j)^{\delta_{jit}}\nonumber\\[-8pt]\\[-8pt]
&&\hspace*{29.8pt}{}\times \{1-p_{j}(\accentset{\leftharpoonup}{X}_{jit,d},
\accentset{\leftharpoonup}{Y}_{jit,d};\theta_{j})\}^{(1-\delta_{jit})}
\bigr]^{\delta_{jit,d}}.\nonumber
\end{eqnarray}
Under some regular conditions, the binary maximum likelihood
estimator $\hat{\theta}_j$ is $\sqrt{n}$-consistent estimator of
$\theta_{j0}$; see \citet{r3} for results on a
related situation. Some guidelines on how to choose models for the
missing propensity are given in Section~\ref{sec8} in the context of the
empirical study. The robustness of the ANOVA tests with respect to
the missing propensity model are discussed in
Sections~\ref{sec3} and~\ref{sec4}.

\section{ANOVA test for covariate effects}\label{sec3}

We consider testing for
$H_{0a}\dvtx\beta_{10}=\beta_{20}=\cdots=\beta_{k0}$ with respect to
the covariates. Let $\pi_{jim}(\theta_j) =\prod_{l=m-d}^{m}
p_j(\accentset{\leftharpoonup}{X}_{jil,d}$, $\accentset{\leftharpoonup
}{Y}_{jil,d};\theta_j)$ be the overall missing propensity for the
$(j,i)$th subject up to time $t_{jim}$. To remove the nonparametric
part in (\ref{eq:model}), we first estimate the nonparametric function
$g_{j0}(t)$. If $\beta_{j0}$ and $\gamma_{j0}$ were known, $g_{j0}(t)$
would be estimated by
%
%
\begin{equation}\label{eq:g0est}
\hat{g}_j(t;\beta_{j0})=\sum_{i=1}^{n_{j}}\sum
_{m=1}^{T_j}w_{jim,h}(t)(Y_{jim}-X_{jim}^{\tau}\beta
_{j0}-M_{jim}^\tau\gamma_{j0}),
\end{equation}
where
%
%
\begin{equation}
\label{eq3.2}
w_{jim,h_j}(t)=\frac{(\delta_{jim}/\pi_{jim}(\hat{\theta
}_j))K_{h_j}(t_{jim}-t)}{\sum_{s=1}^{n_{j}}\sum_{l=1}^{T_j}(\delta
_{jsl}/\pi_{jsl}(\hat{\theta}_j))K_{h_j}(t_{jsl}-t)}
\end{equation}
is a kernel weight that has been inversely weighted by the
propensity $\pi_{jim}(\hat{\theta}_{j})$ to correct for selection
bias due to the missing values. In (\ref{eq3.2}), $K$ is a
univariate kernel function which is a symmetric probability
density, $K_{h_j}(t)=K(t/h_j)/h_j$ and $h_j$ is a smoothing
bandwidth. The conventional kernel estimation of $g_{j0}(t)$
without weighting by $\pi_{jsl}(\hat{\theta}_j)$ may be
inconsistent if the missingness depends on the responses $Y_{jil}$,
which can be the case for missing covariates.

Let $A_{jim}$ denote any of $X_{jim},Y_{jim}$ and $M_{jim}$ and define
%
%
\begin{equation}\label{eq:center1}
\tilde{A}_{jim}=A_{jim}-\sum_{i_1=1}^{n_{j}}\sum
_{m_1=1}^{T_j}w_{ji_1m_1,h_j}(t_{jim})A_{ji_1m_1}
\end{equation}
to be the centering of $A_{jim}$ by the kernel conditional mean estimate,
as is commonly exercised in the partially linear regression
[\citet{r11}]. An estimating function for
the $(j,i)$th subject is
\[
Z_{ji}(\beta_j)=\sum_{m=1}^{T_j}\frac{\delta_{jim}}{\pi_{jim}(\hat
{\theta}_j)}\tilde{X}_{jim}(\tilde{Y}_{jim}-\tilde{X}_{jim}^\tau
\beta_j-\tilde{M}_{jim}^\tau\tilde{\gamma}_j),
\]
where $\tilde{\gamma}_j$ is the solution of
\[
\sum_{i=1}^{n_j}\sum_{m=1}^{T_j}\frac{\delta_{jim}}{\pi_{jim}(\hat
{\theta}_j)}\tilde{M}_{jim}(\tilde{Y}_{jim}-\tilde{X}_{jim}^\tau
\beta_{j0}-\tilde{M}_{jim}^\tau\tilde{\gamma}_j)=0
\]
at the true $\beta_{j0}$. Note that 
$E\{Z_{ji}(\beta_{j0})\}=o(1)$. Although it is not exactly zero,
$Z_{ji}(\beta_{j0})$ can still be used as an approximate zero mean
estimating function to formulate an empirical likelihood for
$\beta_{j}$ as follows.

Let $\{p_{ji}\}_{i=1}^{n_j}$ be nonnegative weights allocated to
$\{(X^\tau_{ji},Y_{ji})\}_{i=1}^{n_j}$. The empirical
likelihood for $\beta_j$
is
%
%
\begin{equation}\label{eq:3.2a}
L_{n_j}(\beta_j)=\max\Biggl\{\prod_{i=1}^{n_j}p_{ji}\Biggr\},
\end{equation}
subject to $ \sum_{i=1}^{n_j}p_{ji}=1$ and $
\sum_{i=1}^{n_j}p_{ji}Z_{ji}(\beta_j)=0. $

By introducing a Lagrange multiplier $\lambda_j$ to solve the
above optimization problem and following the standard derivation
in empirical likelihood [\citet{r27}], it can be shown that
%
%
\begin{equation}\label{eq:eq2}
L_{n_j}(\beta_j)=\prod_{i=1}^{n_j} \biggl\{\frac{1}{n_j}\frac
{1}{1+\lambda_j^{\tau}Z_{ji}(\beta_j)} \biggr\},
\end{equation}
where $\lambda_j$ satisfies
%
%
\begin{equation}\label{eq3}
\sum_{i=1}^{n_j}\frac{Z_{ji}(\beta_j)}{1+\lambda_j^{\tau
}Z_{ji}(\beta_j)}=0.
\end{equation}
The maximum\vspace*{-2pt} of $L_{n_j}(\beta_j)$ is
$\prod_{i=1}^{n_j}\frac{1}{n_j}$, achieved at
$\beta_j=\hat{\beta_j}$ and $\lambda_j=0$, where $\hat{\beta_j}$
solves $\sum_{i=1}^{n_j}Z_{ji}(\hat{\beta_j})=0$.

Let $n=\sum_{i=1}^k n_j, n_j/n\rightarrow\rho_j$ for some nonzero
$\rho_j$ as $n\rightarrow\infty$ such that $\sum_{i=1}^k\rho_j=1$.
As the $k$ samples are independent, the joint empirical likelihood
for $(\beta_1,\beta_2,\ldots,\beta_k)$ is
\[
L_n(\beta_1,\beta_2,\ldots,\beta_k)=\prod_{j=1}^k
L_{n_j}(\beta_j).
\]
%
The log likelihood ratio statistic for $H_{0a}$ is
%
%
\begin{eqnarray}
\label{eq:eq_para}
\ell_n
:\!&=&-2 \max_{\beta}\log
L_n(\beta,\beta,\ldots,\beta)+\sum_{j=1}^k n_j\log n_j
\nonumber\\[-8pt]\\[-8pt]
&=& 2
\min_{\beta}\sum_{j=1}^k\sum_{i=1}^{n_j}\log\{1+\lambda_j^{\tau
}Z_{ji}(\beta)\}.\nonumber
\end{eqnarray}

Using a Taylor expansion and the Lagrange multiplier to carry out
the minimization in (\ref{eq:eq_para}),
%
the optimal solution to $\beta$ is
%
%
\begin{equation}
\label{eq6}
\Biggl(\sum_{j=1}^k\Omega_{x_j}B_j^{-1}\Omega_{x_j}
\Biggr)^{-1} \Biggl(\sum_{j=1}^k\Omega_{x_j}B_j^{-1}\Omega_{x_jy_j} \Biggr)
+ o_p(1),
\end{equation}
where $B_j=\lim_{n_j\rightarrow\infty}{(n_jT_j)}^{-1}
\sum_{i=1}^{n_j}E\{Z_{ji}(\beta_{j0})Z_{ji}(\beta_{j0})^{\tau}\}$,
\[
\Omega_{x_{j}}=\frac{1}{\sqrt{n_jT_j}}\sum_{i=1}^{n_j}\sum
_{m=1}^{T_j}E \biggl\{\frac{\delta_{jim}}{\pi_{jim}(\hat{\theta
}_j)}\tilde{X}_{jim}\tilde{X}_{jim}^{\tau} \biggr\}
\]
and
%
\[
\Omega_{x_{j}y_j}=\frac{1}{\sqrt{n_jT_j}}\sum_{i=1}^{n_j}\sum_{m=1}^{T_j}
\frac{\delta_{jim}}{\pi_{jim}(\hat{\theta}_j)}\tilde
{X}_{jim}(\tilde{Y}_{jim}-M_{jim}^\tau\tilde{\gamma}_j).
\]


%



The ANOVA test statistic (\ref{eq:eq_para}) can be viewed as a
nonparametric counterpart of the conventional parametric
likelihood ratio ANOVA test statistic, for instance that
considered in \citet{r10}. Like its parametric
counterpart, the Wilks theorem is maintained for $\ell_n$.
\begin{theorem}\label{th1} If conditions \textup{A1--A4} given in the
\hyperref[app]{Appendix}
hold, then under $H_{0a}$, 
$\ell_n\stackrel{d}{\rightarrow} \chi_{(k-1)p}^2$ as $n \to
\infty$.
\end{theorem}

The theorem
suggests an empirical likelihood ANOVA test that rejects $H_{0a}$ if
$\ell_n >
\chi^2_{(k-1)p, \alpha} $ where $\alpha$ is the significant level
and $\chi^2_{(k-1)p, \alpha}$ is the upper $\alpha$
quantile of the $\chi^2_{(k-1)p}$ distribution.

We next evaluate the power of the empirical likelihood ANOVA test
under a series
of local alternative hypotheses:
%
\[
H_{1a}\dvtx\beta_{j0}=\beta_{10}+c_n n_j^{-1/2} \qquad\mbox{for }
2 \le j \le k,
\]
where $\{ c_n\}$ is a sequence of bounded constants.
Define
$\Delta_\beta=(\beta_{10}^\tau-\beta_{20}^\tau,\beta_{10}^\tau
-\beta_{30}^\tau,\ldots,\beta_{10}^\tau-\beta_{k0}^\tau)^\tau$,
$D_{1j}=\Omega_{x_1}^{-1}\Omega_{x_1y_1}-\Omega_{x_j}^{-1}\Omega_{x_jy_j}$
for $2 \leq j\leq k$ and
$D=(D_{12}^\tau,D_{13}^\tau,\ldots,D_{1k}^\tau)^\tau$. Let
$\Sigma_{D}=\operatorname{Var}(D)$ and
$\gamma^2=\Delta_\beta^\tau\Sigma_{D}^{-1}\Delta_\beta$.
Theorem~\ref{th2} gives the asymptotic distribution of $\ell_n$ under the
local alternatives.
\begin{theorem}\label{th2} Suppose conditions \textup{A1--A4} in the
\hyperref[app]{Appendix} hold, then under $H_{1a}$, $\ell_n\stackrel
{d}{\rightarrow}
\chi_{(k-1)p}^2(\gamma^2)$ as $n \to\infty$.
\end{theorem}

It can be shown that
%
%
\begin{equation}
\label{SigmaD} \Sigma_D = \Omega_{x_1}^{-1} B_1 \Omega_{x_1}^{-1}
{\mathbf1}_{(k-1)} \otimes{\mathbf1}_{(k-1)} +
\operatorname{diag}\{\Omega_{x_2}^{-1}
B_2 \Omega_{x_2}^{-1},\ldots, \Omega_{x_k}^{-1} B_k
\Omega_{x_k}^{-1} \}.\hspace*{-35pt}
\end{equation}
As each $\Omega^{-1}_{x_ j}$ is $O(n^{1/2})$, the noncentral
component $\gamma^2$ is nonzero and bounded. The power of the
$\alpha$ level empirical likelihood ANOVA test is
\[
\beta(\gamma) = P\bigl\{ \chi_{(k-1)p}^2(\gamma^2)>\chi^2_{(k-1)p,\alpha}\bigr\}.
\]
This indicates that the test is able to detect local
departures of size $O(n^{-1/2})$ from~$H_{0a}$, which is the best
rate we can achieve under the local alternative set-up. This is
attained despite the fact that
nonparametric kernel estimation is involved in the formulation,
which has a slower rate of convergence than $\sqrt{n}$, as the
centering in (\ref{eq:center1}) essentially eliminates the effects
of the nonparametric estimation.
\begin{remark}\label{remark1}
When there\vspace*{-1pt} is no missing values, namely all $\delta
_{jim}=1$, we will assign all $\pi_{jim}(\hat{\theta}_j)=1$ and
there is no need to estimate each $\theta_{j0}$. In this case, Theorems
\ref{th1} and~\ref{th2} remain valid.
It is a different matter for estimation as estimation efficiency with
missing values will be less than that without missing values.
\end{remark}
\begin{remark}\label{remark2}
The above ANOVA test is robust against
misspecifying the missing propensity $p_j(\cdot; \theta_{j0})$
provided the missingness\vspace*{-1pt} does not depend on the responses
$\accentset{\leftharpoonup}{Y}_{jit,d}$. This is because despite
the mispecification, the mean of $Z_{ji}(\beta)$ is still
approximately zero and the empirical likelihood formulation remains
valid, as well
as Theorems~\ref{th1} and~\ref{th2}. However, if the missingness
depends on the
responses and if the model is misspecified, Theorems~\ref{th1} and~\ref{th2} will
be affected.
\end{remark}
\begin{remark}\label{remark3}
The empirical likelihood test can be readily modified
for ANOVA
testing on pure
parametric regressions with some parametric time effects
$g_{j0}(t;\eta_j)$ with parameters $\eta_j$.
When there is absence of interaction, we may formulate the
empirical likelihood for $(\beta_j,\eta_j) \in R^{p+s}$ using
\begin{eqnarray*}
Z_{ji}(\beta_j;\eta_j)&=&\sum_{m=1}^{T_j}\frac{\delta_{jim}}{\pi
_{jim}(\hat{\theta}_j)} \biggl(X_{jim}^\tau,\frac{\partial
g^\tau_j(t_{jim};\eta_j)}{\partial\eta_j} \biggr)^\tau\\
&&\hspace*{16.4pt}{}\times
\{Y_{jim}-X_{jim}^\tau\beta_j-g_{j0}(t_{jim};\eta_j) \}
\end{eqnarray*}
as the estimating function for the $(j,i)$th subject. The ANOVA
test can be formulated following the same procedures from
(\ref{eq:eq2}) to (\ref{eq:eq_para}), and both
Theorems~\ref{th1} and~\ref{th2}
remaining valid after updating $p$ with $p+s$ where $s$ is the
dimension of $\eta_j$.\vadjust{\goodbreak}

In our formulation for the ANOVA test here and in the next
section, we rely on the Nadaraya--Watson type kernel estimator. The
local linear kernel estimator may be employed when the boundary
bias may be an issue. However, as we are interested in ANOVA tests
instead of estimation, the boundary bias does not have a leading
order effect.
\end{remark}

\section{ANOVA test for time effects}\label{sec4}
In this section, we consider the ANOVA test for the nonparametric
part
\[
H_{0b}\dvtx
g_{10}(\cdot)=\cdots=g_{k0}(\cdot).
\]
%
We will first formulate an empirical likelihood for $g_{j 0}(t)$ at
each $t$, which
then lead to an overall likelihood ratio for $H_{0b}$. 
We need an estimator of $g_{j 0}(t)$ that is less biased than the
one in (\ref{eq:g0est}). Recall the notation defined in Section~\ref{sec2}:
$\mathbb{X}_{jim}^\tau=(X_{jim}^\tau,M_{jim}^\tau)$ and
$\xi_j^\tau=(\beta_j^\tau,\gamma_j^\tau)$. Plugging-in the
estimator $\hat{\xi}_j$ to (\ref{eq:g0est}), we have
%
%
\begin{equation}\label{eq4.8}
\tilde{g}_j(t)=\sum_{i=1}^{n_j}\sum
_{m=1}^{T_j}w_{jim,h_j}(t)(Y_{jim}-\mathbb{X}_{jim}^{\tau}\hat{\xi
}_{j}).
\end{equation}
It follows that, for any $t\in[0,1]$,
%
%
\begin{eqnarray}
\label{eq4.9}
\tilde{g}_j(t)-g_{j0}(t)&=&\sum_{i=1}^{n_j}\sum_{m=1}^{T_j}w_{jim,h_j}(t)
\{\varepsilon_{ji}(t_{jim})+\mathbb{X}_{jim}^{\tau}
(\xi_j-\hat{\xi}_{j})\nonumber\\[-8pt]\\[-8pt]
&&\hspace*{105.5pt}{}+g_{j0}(t_{jim})-g_{j0}(t) \}.\nonumber
\end{eqnarray}
%
However, there is a bias of order $h_j^2$ in the kernel estimation
since
\[
\sum_{i=1}^{n_j}\sum_{m=1}^{T_j}w_{jim,h_j}(t) \{
g_{j0}(t_{jim})-g_{j0}(t) \}={1\over2}
\biggl\{\int z^2K(z)\,dz\biggr\} g_{j0}''(t)h_j^2+o_p(h_j^2).
\]
If we formulated the empirical likelihood based on $\tilde{g}_j(t)$,
the bias will
contribute to the asymptotic distribution of the ANOVA test
statistic. To avoid that, we use
the bias-correction method proposed in \citet{r45} so that
the estimator of $g_{j0}$ is
\[
\hat{g}_{j}(t)=\sum_{i=1}^{n_j}\sum_{m=1}^{T_j}w_{jim,h_j}(t)\bigl\{
Y_{jim}-\mathbb{X}_{jim}^{\tau}\hat{\xi}_{j}-\bigl(\tilde
{g}_{j}(t_{jim})-\tilde{g}_{j}(t)\bigr)\bigr\}.
\]
%
Based on this modified estimator $\hat{g}_{j}(t)$, we define the
auxiliary variable
\begin{eqnarray*}
R_{ji}\{g_j(t)\}&=&\sum_{m=1}^{T_j}\frac{\delta_{jim}}{\pi
_{jim}(\hat{\theta}_j)}
K \biggl(\frac{t_{jim}-t}{h_j} \biggr)\\
&&\hspace*{16.4pt}{} \times
\bigl\{Y_{jim}-\mathbb{X}_{jim}^{\tau}\hat{\xi}_{j}-g_j(t)-
\bigl(\tilde{g}_j(t_{jim})-\tilde{g}_j(t) \bigr) \bigr\}
\end{eqnarray*}
for empirical likelihood formulation. At true function $g_{j0}(t)$,
$E[R_{ji}\{g_{j0}(t)\} ]=o(1)$.

Using a similar procedure to $L_{n_j}(\beta_j)$ as given in
(\ref{eq:eq2}) and (\ref{eq3}),
the empirical likelihood for $g_{j0}(t)$ is
\[
L_{n_j}\{g_{j0}(t)\}=\max\Biggl\{\prod_{i=1}^{n_j}p_{ji} \Biggr\}
\]
subject to $ \sum_{i=1}^{n_j} p_{ji}=1$ and $
\sum_{i=1}^{n_j}p_{ji} R_{ji}\{g_j(t)\}=0$. The latter is obtained
in a similar fashion as we obtain (\ref{eq:eq2}) by introducing
Lagrange multipliers so that
\[
L_{n_j} \{ g_{j 0}(t)\} =
\prod_{i=1}^{n_j} \biggl\{\frac{1}{n_j}\frac{1}{1+\eta_j(t) R_{ji}
\{g_{j 0}(t)\}} \biggr\},
\]
where $\eta_j(t)$ is a Lagrange multiplier that satisfies
%
%
\begin{equation}\label{2.13}
\sum_{i=1}^{n_j}\frac{ R_{ji} \{g_{j 0}(t)\} }{1+\eta_j(t) R_{ji}
\{g_{j 0}(t)\} }=0.
\end{equation}
The log empirical likelihood ratio for
$g_{10}(t)=\cdots=g_{k0}(t):=g(t)$, say, is
%
%
\begin{equation}
\label{eq:eq_nonpara} {\mathcal{L}}_n(t)=2
\min_{g(t)}\sum_{j=1}^k\sum_{i=1}^{n_j}\log\bigl(1+\eta_j(t)R_{ji}\{
g(t)\}\bigr),
\end{equation}
which is analogues of $\ell_n$ in (\ref{eq:eq_para}). As shown in
the proof of Theorem~\ref{th3} given in the supplement article
[\citet{r4}], the leading order term of the ${\mathcal{L}}_n(t)$
is a
studentized version of the distance
\[
\bigl( \hat{g}_1(t) - \hat{g}_2(t), \hat{g}_1(t) - \hat{g}_3(t),\ldots, \hat
{g}_1(t) - \hat{g}_k(t) \bigr),
\]
namely between $\hat{g}_1(t)$ and the other $\hat{g}_j(t) (j\neq
1)$.
This motivates us to propose using
%
%
\begin{equation}
{\mathcal{T}}_n=\int_{0}^1{\mathcal{L}}_n(t)\varpi(t)\,dt
\end{equation}
to test for the equivalence of $\{ g_{j 0}(\cdot)\}_{j=1}^k$,
where $\varpi(t)$ is a probability weight function over $[0,1]$.

To define the asymptotic distribution of ${\mathcal{T}}_n$, we assume
without loss of generality that for each
$h_j$ and $T_j$, $j=1,\ldots,k$, there exist fixed finite
positive constants $\alpha_j$ and $b_j$ such that $\alpha_j T_j=
T$ and $b_j h_j = h$ for some $T$ and $h$ as $h \to0$.
Effectively, $T$ is the smallest common multiple of
$T_1,\ldots,T_k$. Let $K^{(2)}_c(t)=\int K(w)K(t-cw)\,dt$ and
$K^{(4)}_c(0)=\int
K^{(2)}_{c}(w\sqrt{c})K^{(2)}_{1/c}(w/\sqrt{c})\,dw$. For
$c=1$,\vadjust{\goodbreak}
we\vspace*{-1pt}
resort to the standard notations of $K^{(2)}(t)$ and $K^{(4)}(0)$
for $K^{(2)}_1(t)$ and $K^{(4)}_1(0)$, respectively. For each
treatment $j$, let $f_j$ be the super-population density of the
design points $\{t_{jim}\}$. Let $a_j=\rho_j^{-1}\alpha_j$,
\[
W_j(t)=\frac{f_j(t)/\{a_jb_j\sigma_{\varepsilon
j}^2\}}{\sum_{l=1}^k f_l(t)/\{a_lb_l\sigma_{\varepsilon l}^2\}}
\]
and $V_j(t)=K^{(2)}(0)\sigma_{\varepsilon j}^2f_j(t)$ where
$\sigma_{\varepsilon
j}^2=\frac{1}{n_jT_j}\sum_{i=1}^{n_j}\sum_{m=1}^{T_j}E\{\frac
{\varepsilon^2_{jim}}{\pi_{jim}(\theta_{j0})}\}$.
Furthermore, we define
\begin{eqnarray*}
\Lambda(t)&=&\sum_{j=1}^kb_j^{-1}K^{(4)}(0)\bigl(1-W_j(t)\bigr)^2\\
&&{}+\sum_{j\neq
j_1}^k(b_jb_{j_1})^{-
{1/2}}K^{(4)}_{b_j/{b_{j_1}}}(0)W_j(t)W_{j_1}(t)
\end{eqnarray*}
and
\begin{eqnarray*}
\mu_1&=&\int_{0}^1 \Biggl[\sum_{j=1}^kb_j^{-
{1/2}}V_j^{-1}(t)f_j^2(t)\Delta_{nj}^{2}(t)\\
&&\hspace*{18.2pt}{} - \Biggl(\sum_{s=1}^kb_s^{-{1/4}}V_{s}^{-
{1/2}}(t)W_s^{{1/2}}(t)f_s(t)\Delta_{ns}(t) \Biggr)^2
\Biggr]\varpi(t)\,dt.
\end{eqnarray*}

We consider a sequence of local alternative hypotheses:
%
%
\begin{equation}\label{eq:localg}
g_{j0}(t)=g_{10}(t)+C_{jn}\Delta_{nj}(t),
\end{equation}
where $C_{jn}=(n_jT_j)^{-1/2}h_j^{-1/4}$ for $j= 2,\ldots, k$ and
$\{\Delta_{nj}(t)\}_{n \ge1}$ is a sequence of uniformly bounded
functions.
%
%
\begin{theorem}\label{th3} Assume conditions \textup{A1--A4} in the
\hyperref[app]{Appendix} and
$h=O(n^{-1/5})$, then under (\ref{eq:localg}),
\[
h^{-1/2}({\mathcal{T}}_n-\mu_0)\stackrel{d}{\rightarrow}
N(0,\sigma^2_0),
\]
where $\mu_0=(k-1)+h^{1/2}\mu_1$ and
$\sigma_0^2=2K^{(2)}(0)^{-2}\int_0^1\Lambda(t)\varpi^2(t)\,dt$.
\end{theorem}

We note that under $H_{0b}\dvtx g_{10}(\cdot)=\cdots=g_{k0}(\cdot)$,
$\Delta_{nj}(t)=0$ which yields $\mu_1=0$ and
%
\[
h^{-1/2} \{ {\mathcal{T}}_n- (k-1)\} \stackrel{d}{\rightarrow}
N(0,\sigma^2_0).
\]
This may lead to an asymptotic test at a nominal significance
level $\alpha$ that rejects $H_{0b}$ if
%
%
\begin{equation}\label{eq:asytest}
{\mathcal{T}}_n \ge
h^{1/2} \hat{\sigma}_0 z_{\alpha} + (k-1),\vadjust{\goodbreak}
\end{equation}
where $z_{\alpha}$ is the upper $\alpha$ quantile of $N(0,1)$ and
$\hat{\sigma}_0$ is a consistent estimator of $\sigma_0$. The
asymptotic power of the test under the local alternatives
is $1-\Phi(z_\alpha-\frac{\mu_1}{\sigma_0})$, where
$\Phi(\cdot)$ is the standard normal distribution function. This
indicates that the test is powerful in differentiating null
hypothesis and its local alternative at the convergence rate
$O(n_j^{-1/2} h_j^{-1/4})$ for $C_{jn}$. The rate is the best when
a single bandwidth is used [\citet{r12}].

If all the $h_j (j=1,\ldots,k)$ are the same, 
the asymptotic variance $\sigma_0^2=2(k-1)K^{(2)}(0)^{-2}
K^{(4)}(0)\int_0^1\varpi^2(t)\,dt$, which means that the test
statistic under $H_{0b}$ is asymptotic pivotal.
However, when the bandwidths are not the
same, which is most likely as different treatments may
require different amount of smoothness in the estimation of
$g_{j0}(\cdot)$, the asymptotical pivotalness of ${\mathcal{T}}_n$ is
no longer available, and estimation of $\sigma_0^2$ is needed for
conducting the asymptotic test in (\ref{eq:asytest}). We will
propose a test based on a bootstrap calibration to the distribution
of ${\mathcal{T}}_n$ in Section~\ref{sec6}.
\begin{remark}\label{remark4}
Similar to Remarks~\ref{remark1} and~\ref{remark2} made on the ANOVA
tests for the covariate effects, the proposed ANOVA test for the
nonparametric baseline functions (Theorem~\ref{th3}) remains valid in the
absence of missing values or if the missing propensity is
misspecified as long as the responses do not contribute to the
missingness.
\end{remark}
\begin{remark}\label{remark5}
We note that the proposed test is not affected by
the within-subject dependent structure (the longitudinal aspect)
due to the fact that the formulation of the empirical likelihood is
made for each
subject. This is clearly shown in the construction of
$R_{ji}\{g_j(t)\}$ and by the fact that the nonparametric
functions can be separated from the covariate effects in the
semiparametric model. Again this would be changed if we are
interested in estimation as the correlation structure in the
longitudinal data will affect the estimation efficiency. However,
the test will be dependent on the choice of the weight function
$\varpi(\cdot)$, and $\{\alpha_j\}$, $\{\rho_j\}$ and $\{ b_j\}$,
the relative ratios among $\{T_j\}$, $\{n_j\}$ and $\{ h_j\}$.
\end{remark}
\begin{remark}\label{remark6}
The ANOVA test statistics for the time effects for the
semiparametric model can be readily modified to obtain ANOVA test for purely
nonparametric regression by simply setting $\hat{\xi}_j=0$ in
the formulation of the test statistic ${\mathcal{L}}_n(t)$. In this
case, the model (\ref{eq:model}) takes the form
%
%
\begin{equation}
\label{eq:model_nonpara}
Y_{ji}(t)=g_j(X_{ji}(t),t)+\varepsilon_{ji}(t),
\end{equation}
where $g_j(\cdot)$ is the unknown nonparametric function of
$X_{ji}(t)$ and $t$. The proposed ANOVA test can be viewed as
generalization of the tests considered in \citet{r25},
\citet{r29}
and \citet{r38} by considering both the
longitudinal and missing aspects. See also
\citet{r0} for a two sample test for the equivalence of two
probability densities.
\end{remark}

\section{Tests on interactions}\label{sec5}

Model (\ref{eq:model0}) contains an interactive term\break $M(X_{jim}, t)$
that is flexible in prescribing the interact between $X_{jim}$ and the
time, as long as the positive definite condition in condition A3 is
satisfied. In this section, we propose tests for the presence of the
interaction in the $j$th treatment and the ANOVA hypothesis on the
equivalence of the interactions among the treatments.

We firstly consider testing $H_{0c}\dvtx\gamma_{j0}=0$
vs. $H_{1c}\dvtx\gamma_{j0}\neq0$ for a fixed $j$. In the
formulation of the empirical likelihood for $\gamma_{j0}$, we treat
$M_{jim}=M(X_{jim}, t)$ as a covariates with the same role like
$X_{jim}$ in the previous section when we constructed empirical
likelihood for
$\beta_{j 0}$. For this purpose, we define estimating equations
for $\gamma_{j0}$
%
%
\begin{equation}
\phi_{ji}(\gamma_{j0})=\sum_{m=1}^{T_j}\frac{\delta_{jim}}{\pi
_{jim}(\hat{\theta}_j)}
\tilde{M}_{jim}(\tilde{Y}_{jim}-\tilde{X}_{jim}^\tau\tilde{\beta
}_j-\tilde{M}_{jim}^\tau\gamma_{j0}),
\end{equation}
where
%
%
\begin{eqnarray}
\tilde{\beta}_{j}&=& \Biggl\{\sum_{i=1}^{n_j}\sum_{m=1}^{T_j}\frac
{\delta_{jim}}{\pi_{jim}(\hat{\theta}_j)}
\tilde{X}_{jim}\tilde{X}_{jim}^\tau\Biggr\}^{-1}
\nonumber\\[-8pt]\\[-8pt]
&&{}\times\sum_{i=1}^{n_j}\sum_{m=1}^{T_j}\frac{\delta_{jim}}{\pi_{jim}(\hat
{\theta}_j)}\tilde{X}_{jim}
(\tilde{Y}_{jim}-\tilde{M}_{jim}^\tau\gamma_{j0})\nonumber
\end{eqnarray}
is the
``estimator'' of $\beta_j$ at the true $\gamma_{j0}$. Similar to
establishing $\ell_{n_j}(\beta_j)$ in Section~\ref{sec3}, the log-empirical
likelihood for $\gamma_{j0}$ can be written as
\[
\ell^\gamma_{n_j}(\gamma_j)=2\sum_{i=1}^{n_j}\log\{1+\Lambda
_j'\phi_{ji}(\gamma_{j})\},
\]
where the Lagrange multipliers $\Lambda_j$ satisfies
%
%
\begin{equation}\label{eq6:A1}\sum_{i=1}^{n_j}\frac{\phi
_{ji}(\gamma
_{j})}{1+\Lambda_j'\phi_{ji}(\gamma_{j})}=0.
\end{equation}
To test for $H_{0d}\dvtx\gamma_{10}=\gamma_{20}=\cdots=\gamma_{k0}$
vs. $H_{1d}\dvtx\gamma_{i0}\neq\gamma_{j0}$ for some
$i\neq j$,
we construct the joint empirical likelihood ratio
%
%
\begin{equation}\label{eq:eq_interact}\ell^\gamma_n := 2
\min_{\gamma}\sum_{j=1}^k\sum_{i=1}^{n_j}\log\{1+\Lambda_j^{\tau
}\phi_{ji}(\gamma)\},
\end{equation}
where $\Lambda_j$ satisfy (\ref{eq6:A1}).

The asymptotic distributions of the empirical likelihood ratios
$\ell^\gamma_{n_j}({0})$ and $\ell^\gamma_n$ under the null
hypotheses are given in the next theorem whose proofs will not be
given as they follow the same routes in the proof of Theorem~\ref{th1}
by exchanging $X_{jim}$ and $\beta_{j0}$ with $M_{jim}$ and
$\gamma_{j0}$, respectively.
\begin{theorem}\label{th5} Under conditions \textup{A1--A4} given in the
\hyperref[app]{Appendix}, then \textup{(i)}~under $H_{0c}$,
$\ell^\gamma_{n_j}(\mathbf{0})\stackrel{d}{\rightarrow} \chi_{q}^2$
as $n_j \to\infty$; \textup{(ii)} under $H_{0d}$, $ \ell^\gamma_n \stackrel
{d} \to\chi^2_{(k-1) q}$ as $n \to\infty$.
\end{theorem}

Based on Theorem~\ref{th5}, an $\alpha$-level empirical likelihood
ratio test for
the presence of the interaction in the $j$th sample rejects
$H_{0c}$ if $\ell^\gamma_{n_j}(\mathbf{0})>\chi_{q,\alpha}^2$, and
the ANOVA test for the equivalence of the interactive effects
rejects $H_{0 d}$ if $\ell^\gamma_n > \chi^2_{(k-1) q, \alpha}$.
The ANOVA test for $H_{0d}$ has a similar local power performance
as that described after Theorem~\ref{th2} for the ANOVA test regarding
$\beta_{j0}$ in Section~\ref{sec3}. The power properties of the test for
$H_{0c}$ can be established using a much easier method.

We have assumed parametric models for the interaction in model
(\ref{eq:model0}). A semiparametric model would be employed to
model the interaction given that the model for the time effect is
nonparametric. The parametric interaction is a simplification and
avoids some of the involved technicalities associated with a
semiparametric model.





\section{Bootstrap calibration}\label{sec6}

To avoid direct estimation of $\sigma_0^2$ in Theorem~\ref{th3} and to
speed up the convergence of ${\mathcal{T}}_n$, we resort to the
bootstrap. While the wild bootstrap [\citet{r42}, \citet{r21} and
\citet{r12}] originally proposed for parametric
regression without missing values has been modified by \citet{r35}
to take into account missing values, we extend it
further to suit the longitudinal feature.

Let $\vec{t}_j^o$ and $\vec{t}_j^m$ be the sets of the time points
with full and missing observations, respectively. According to
model (\ref{Secondmodel}), we impute a missing $X_{ji}(t)$ from
$\hat{X}_{ji}(t), t\in\vec{t}_j^o$, so that for any
$t\in\vec{t}_j^m$
%
%
\begin{equation}
\label{eq4.15}
\hat{X}_{ji}(t)=\sum_{i=1}^{n_j}\sum_{m=1}^{T_j}w_{jim,h_j}(t)X_{jim},
\end{equation}
where $w_{jim,h_j}(t)$ is the kernel weight defined in
(\ref{eq3.2}).

To mimic the heteroscedastic and correlation structure in the
longitudinal data, we estimate the covariance matrix for each
subject in each treatment. Let
\[
\hat{\varepsilon}_{jim}=Y_{jim}-\mathbb{X}_{jim}^{\tau}\hat{\xi
}_{j}-\hat{g}_{j}(t_{jim}).
\]
%
An estimator of $\sigma_j^2(t)$, the variance of
$\varepsilon_{ji}(t)$, is
$\hat{\sigma}^2_j(t)=\sum_{i=1}^{n_j}\sum
_{m=1}^{T_j}w_{jim,h_j}(t)\times\hat{\varepsilon}_{jim}^2$
and an estimator of $\rho_j(s,t)$, the correlation coefficient
between $\varepsilon_{ji}(t)$ and $\varepsilon_{ji}(s)$ for $s\neq
t$, is
\[
\hat{\rho}_j(s,t)=\sum_{i=1}^{n_j}\sum_{m\neq
m'}^{T_j}H_{jim,m'}(s,t)\hat{e}_{jim}\hat{e}_{jim'},
\]
where
$\hat{e}_{jim}=\hat{\varepsilon}_{jim}/\hat{\sigma}_j(t_{jim})$,
\[
\hspace*{-5pt}H_{jim,m'}(s,t)=\frac{\delta_{jim}\delta
_{jim'}K_{b_j}(s-t_{jim})K_{b_j}(t-t_{jim'})/\pi_{jim,m'}(\hat{\theta}_j)}
{\sum_{i=1}^{n_j}\sum_{m\neq
m'}\delta_{jim}\delta_{jim'}K_{b_j}(s-t_{jim})K_{b_j}(t-t_{jim'})/\pi
_{jim,m'}(\hat{\theta}_j)}
\]
and
$\pi_{jim,m'}(\hat{\theta}_j)=\pi_{jim}(\hat{\theta}_j)\pi
_{jim'}(\hat{\theta}_j)$
if $|m-m'|>d;$
$\pi_{jim,m'}(\hat{\theta}_j)=\pi_{jim_b}(\hat{\theta}_j)$ if
$|m-m'|\leq d$ where $m_b=\max(m,m')$. Here $b_j$ is a smoothing
bandwidth which may be different from the bandwidth $h_j$ for
calculating the test statistics ${\mathcal{T}}_n$ [\citet{r7}].
Then, the covariance $\Sigma_{ji}$ of
$\varepsilon_{ji}=(\varepsilon_{ji1},\ldots,\varepsilon
_{jiT_j})^\tau$
is estimated by $\hat{\Sigma}_{ji}$ which has $\hat{\sigma
}^2_j(t_{jim})$ as
its $m$th diagonal element and
$\hat{\rho}_j(t_{jik},t_{jil})\hat{\sigma}_j(t_{jik})\hat{\sigma
}_j(t_{jil})$
as its $(k,l)$th element for $k\neq l$.

Let $Y_{ji},\delta_{ji},t_{ji}$ be the vector of random variables
of the $(j,i)$th subject,
$X_{ji}=(X_{ji}(t_{ji1}),\ldots,X_{ji}(t_{jiT_j}))^\tau$ and
$g_{j0}(t_{sl})=(g_{j0}(t_{sl1}),\ldots,
g_{j0}(t_{slT_k}))^\tau$,\break
where $s$ may be different from $j$. Let
$X_{ji}^c=\{X_{ji}^o,\hat{X}_{ji}^m\}$, where $X_{ji}^o$ contains
observed $X_{ji}(t)$ for $t_j\in\vec{t}^o$ and $\hat{X}_{ji}^m$
collects the imputed $X_{ji}(t)$ for $t \in\vec{t}_j^m$
according to~(\ref{eq4.15}). Plugging the value of $X_{ji}^c$, we
get $M_{ji}^c=\{ M_{ji}^o, \hat{M}_{ji}^m\}$, the observed and the
imputed interactions for $(j,i)$th subject and then
$\mathbb{X}_{ji}^c$.

The proposed bootstrap procedure consists of the following steps:

\textit{Step} 1. Generate a bootstrap re-sample
$\{Y^*_{ji},\mathbb{X}_{ji}^c,\delta^*_{ji},t_{ji}\}$ for the
$(j,i)$th subject by
\[
Y^*_{ji}={\mathbb{X}_{ji}^c}^\tau\hat{\xi}_j+\hat
{g}_{1}(t_{ji})+\hat{\Sigma}_{ji}e^*_{ji},
\]
where $e^*_{ji}$'s are i.i.d. random vectors simulated from a
distribution satisfying $E(e^*_{ji})=0$ and
$\operatorname{Var}(e^*_{ji})=I_{T_j}$,
$\delta^*_{jim}\sim\operatorname{Bernoulli}(\pi_{jim}(\hat{\theta}_j))$
where $\hat{\theta}_j$ is estimated based on the original sample
as given in (\ref{BinLikehood}). Here, $\hat{g}_{1}(t_{ji})$ is
used as the common nonparametric time effect to mimic the null
hypothesis $H_{0b}$.

\textit{Step} 2. For each treatment $j$, we reestimate $\xi_j$,
$\theta_j$ and $g_j(t)$ based on the resample
$\{Y^*_{ji},\mathbb{X}_{ji}^c,\delta^*_{ji},t_{ji}\}$ and denote
them as $\hat{\xi}_j^*$, $\hat{\theta}_j^*$ and
$\hat{g}^*_{j}(t)$.
The bootstrap version of $R_{ji}\{{g}_1(t)\}$ is
\begin{eqnarray*}
R_{ji}^*\{\hat{g}_1(t)\}&=&\sum_{m=1}^{T_j}\frac{\delta_{jim}^*}{\pi
_{jim}(\hat{\theta}^*_j)}
K \biggl(\frac{t_{jim}-t}{h_j} \biggr)\\
&&\hspace*{16.8pt}{} \times\bigl\{Y_{jim}^*-\mathbb{X}_{jim}^{\tau}\hat{\xi
}_j^*-\hat{g}_1(t)-\{\hat{g}^*_j(t_{jim})-\hat{g}^*_j(t)\} \bigr\}
\end{eqnarray*}
and use it to substitute $R_{ji}\{g_j(t)\}$ in the formulation of
${\mathcal{L}}_n(t)$, we obtain ${\mathcal{L}}^*_n(t)$ and then
${\mathcal{T}}_n^*=\int{\mathcal{L}}^*_n(t)\varpi(t)\,dt$.

\textit{Step} 3. Repeat the above two steps $B$ times for a large
integer $B$ and obtain $B$ bootstrap values
$\{{\mathcal{T}}_{nb}^*\}_{b=1}^B$. Let $\hat{t}_\alpha$ be the
$1-\alpha$ quantile of $\{{\mathcal{T}}_{nb}^*\}_{b=1}^B$, which is a
bootstrap estimate of the $1-\alpha$ quantile of ${\mathcal{T}}_n$.
Then, we reject the null hypothesis $H_{0b}$ if
${\mathcal{T}}_n>\hat{t}_\alpha$.

The following theorem justifies the bootstrap procedure.
%
%
\begin{theorem}\label{th4} Assume conditions \textup{A1--A4} in the
\hyperref[app]{Appendix}\vspace*{1pt}
hold and $h=O(n^{-1/5})$. Let ${\mathcal{X}}_n$ denote the original
sample, $h$ and $\sigma_0^2$ be defined as in Theorem~\ref{th3}. The
conditional distribution of $h^{-1/2}({\mathcal{T}}_n^*-\mu_0)$ given
${\mathcal{X}}_n$ converges to $N(0,\sigma_0^2)$ almost surely,
namely,
\[
h^{-1/2}\{
{\mathcal{T}}_n^*-(k-1)\}|{\mathcal{X}}_n\stackrel{d}{\rightarrow}
N(0,\sigma^2_0) \qquad\mbox{a.s.}
\]
%
\end{theorem}

\section{Simulation results}\label{sec7}

In this section, we report results from simulation studies which
were designed to confirm the proposed ANOVA tests proposed in the
previous sections. We simulated data from the following
three-treatment model:
%
%
\begin{eqnarray}
\label{eq5.4}
Y_{jim}&=&X_{jim}\beta_{j}+M_{jim}\gamma_j+g_j(t_{jim})+\varepsilon_{jim}
\quad\mbox{and}\nonumber\\[-8pt]\\[-8pt]
X_{jim}&=&2-1.5t_{jim}+u_{jim},\nonumber
\end{eqnarray}
where $M_{jim}=t_{jim}\times(X_{jim}-1.5)^2$,
$\varepsilon_{jim}=e_{ji}+\nu_{jim}$, $u_{jim}\sim
N(0,\sigma_{a_j}^2)$, $e_{ji}\sim N(0, \sigma_{b_j}^2)$ and
$\nu_{jim}\sim N(0, \sigma_{c_j}^2)$ for $j=\{1,2,3\}$,
$i=1,\ldots, n_j$ and $m=1,\ldots,T_j$. This structure used to
generate $\{\varepsilon_{jim}\}_{m=1}^{T_j}$ ensures dependence
among the repeated measurements $\{Y_{jim}\}$ for each subject
$i$. The correlation between $Y_{jim}$ and $Y_{jil}$ for any
$m\neq l$ is $\sigma_{b_j}^2/(\sigma_{b_j}^2+\sigma_{c_j}^2)$. The
time points $\{t_{jim}\}_{m=1}^{T_j}$ were obtained by first
independently generating uniform$[0,1]$ random variables and then
sorted in the ascending order. We set the number of repeated
measures $T_j$ to be the same, say $T$, for all three treatments;
and chose $T=5$ and 10, respectively. The standard deviation
parameters in (\ref{eq5.4}) were $\sigma_{a_1}=0.5,
\sigma_{b_1}=0.5, \sigma_{c_1}=0.2$ for the first treatment,
$\sigma_{a_2}=0.5, \sigma_{b_2}=0.5, \sigma_{c_2}=0.2$ for the
second and $\sigma_{a_3}=0.6, \sigma_{b_3}=0.6, \sigma_{c_3}=0.3$
for the third.

The parameters and the time effects for the three treatments were:
%
\begin{longlist}
\item[\textit{Treatment} 1:] $\beta_1=2, \gamma
_1=1, g_1(t)=2\sin(2\pi t)$; 
\item[\textit{Treatment} 2:] $\beta_2=2+D_{2n}, \gamma
_2=1+D_{2n}, g_2(t)=2\sin(2\pi t)-\Delta_{2n}(t);$
\item[\textit{Treatment} 3:] $\beta_3=2+D_{3n}, \gamma
_3=1+D_{3n}, g_3(t)=2\sin(2\pi t)-\Delta_{3n}(t)$.
\end{longlist}
We designated different values of $D_{2n}, D_{3n},
\Delta_{2n}(t)$
and $\Delta_{3n}(t)$ in the evaluation of the size and the power,
whose details will be reported shortly.

We considered two missing data mechanisms. In the first mechanism
(I), the missing propensity was
%
%
\begin{equation}
\label{eq5.1}
\operatorname{logit}\{P(\delta_{jim}=1|\delta_{jim,m-1}=1,X_{ji},Y_{ji})\}=\theta_j
X_{ji(m-1)} \qquad\mbox{for } m>1,\hspace*{-25pt}
\end{equation}
which is not dependent on the response $Y$, with $\theta_{1}=3,
\theta_{2}=2$ and $\theta_{3}=2 $.
In the second mechanism (II),
%
%
\begin{eqnarray}
&& \operatorname{logit}\{P(\delta_{jim}=1|\delta_{jim,m-1}=1,X_{ji},Y_{ji})\}
\nonumber\\[-8pt]\\[-8pt]
&&\qquad =\cases{
\theta_{j1}X_{ji(m-1)}+\theta_{j2}\bigl\{Y_{ji(m-1)}-Y_{ji(m-2)}\bigr\},&\quad
if $m>2$,\cr
\theta_{j1}X_{ji(m-1)}, &\quad if $m=2$;}\nonumber
\end{eqnarray}
which is influenced by both covariate and response, with
$\theta_1=(\theta_{11},\theta_{12})^\tau=(2,-1)^\tau,
\theta_2=(\theta_{21},\theta_{22})^\tau=(2,-1.5)^\tau$ and
$\theta_3=(\theta_{31},\theta_{32})^\tau=(2,-1.5)^\tau$. In both
mechanisms, the first observation $(m=1)$
for each subject was always observed as we have assumed earlier.

We used the Epanechnikov kernel $K(u)=0.75(1-u^2)_{+}$ throughout
the simulation where $(\cdot)_{+}$ stands for the positive part of
a function. The bandwidths were chosen by the ``leave-one-subject''
out cross-validation. Specifically, we chose the bandwidth $h_j$
that minimized the cross-validation score functions
\[
\sum_{i=1}^{n_j}\sum_{m=1}^{T_j}\frac{\delta_{jim}}{\pi_{jim}(\hat
{\theta}_j)}
\bigl(Y_{jim}-X_{jim}^\tau\hat{\beta}_j^{(-i)}-M_{jim}^\tau\hat{\gamma
}_j^{(-i)}-\hat{g}^{(-i)}_j(t_{jim})\bigr)^2,
\]
where $\hat{\beta}_j^{(-i)}$, $\hat{\gamma}_j^{(-i)}$ and
$\hat{g}^{(-i)}_j(t_{jim})$ were the corresponding estimates
without using observations of the $i$th subject. The
cross-validation was used to choose an optimal bandwidth for
representative data sets and fixed the chosen bandwidths in the
simulations with the same sample size. We fixed the number of
simulations to be 500.

The average missing percentages based on 500 simulations for the
missing mechanism I were 8\%, 15\% and 17\% for treatments 1--3,
respectively, when $T=5$, and were 16\%, 28\% and 31\% when $T=10$.
In the missing mechanism~II,
the average missing percentages
were 10\%, 8\% and 15\% for $T=5$, and 23\%, 20\% and 36\% for
$T=10$, respectively.

%
\begin{table}
\caption{Empirical size and power of the 5\% ANOVA test for
$H_{0a}\dvtx\beta_{10}=\beta_{20}=\beta_{30}$}
\label{table1}
\begin{tabular*}{\tablewidth}{@{\extracolsep{\fill}}lccclcccccc@{}}
\hline
\multicolumn{3}{@{}c}{\textbf{Sample size}} & & & & \multicolumn{2}{c}{\textbf{Missingness}} & &
\multicolumn{2}{c@{}}{\textbf{Missingness}}\\[-4pt]
\multicolumn{3}{@{}c}{\hrulefill} & & & &\multicolumn{2}{c}{\hrulefill} & &
\multicolumn{2}{c@{}}{\hrulefill}\\
$\bolds{n_1}$ & $\bolds{n_2}$ & $\bolds{n_3}$ & $\bolds{D_{2n}}$ & \multicolumn{1}{c}{$\bolds{D_{3n}}$} & $\bolds{T}$
& \textbf{I} & \textbf{II} & $\bolds{T}$ & \textbf{I} & \textbf{II}\\
\hline
\phantom{0}60 & \phantom{0}65 & \phantom{0}55 & 0.0 & 0.0 (size) & 5 & 0.042 & 0.050 & 10 & 0.046 &
0.044 \\
& & & 0.2 & 0.0 & & 0.192 & 0.254 & & 0.408 & 0.434 \\
& & & 0.3 & 0.0 & & 0.548 & 0.630 & & 0.810 & 0.864 \\
& & & 0.0 & 0.2 & & 0.236 & 0.214 & & 0.344 & 0.354 \\
& & & 0.0 & 0.3 & & 0.508 & 0.546 & & 0.714 & 0.722 \\
& & & 0.2 & 0.2 & & 0.208 & 0.262 & & 0.446 & 0.458 \\
& & & 0.2 & 0.3 & & 0.412 & 0.440 & & 0.680 & 0.698 \\
& & & 0.3 & 0.2 & & 0.426 & 0.490 & & 0.728 & 0.728 \\
& & & 0.3 & 0.3 & & 0.594 & 0.620 & & 0.836 & 0.818 \\
100 & 110 & 105 & 0.0 & 0.0 (size) & 5 & 0.052 & 0.054 & 10 & 0.042 &
0.038 \\
& & & 0.2 & 0.0 & & 0.426 & 0.470 & & 0.686 & 0.718 \\
& & & 0.3 & 0.0 & & 0.854 & 0.854 & & 0.964 & 0.974 \\
& & & 0.0 & 0.2 & & 0.406 & 0.444 & & 0.612 & 0.568 \\
& & & 0.0 & 0.3 & & 0.816 & 0.836 & & 0.936 & 0.910 \\
& & & 0.2 & 0.2 & & 0.404 & 0.480 & & 0.674 & 0.686 \\
& & & 0.2 & 0.3 & & 0.744 & 0.694 & & 0.944 & 0.882 \\
& & & 0.3 & 0.2 & & 0.712 & 0.768 & & 0.922 & 0.920 \\
& & & 0.3 & 0.3 & & 0.824 & 0.814 & & 0.972 & 0.970 \\
\hline
\end{tabular*}   \vspace*{-2pt}
\end{table}

For the ANOVA test for $H_{0a}\dvtx\beta_{10}=\beta_{20}=\beta_{30}$
with respect to the covariate effects, three values of $D_{2n}$
and $D_{3n}$: 0, 0.2 and 0.3, were used, respectively, while
$\Delta_{2n}(t)=\Delta_{3n}(t)=0$. Table~\ref{table1}
summarizes the empirical size and power of the proposed EL ANOVA
test with 5\% nominal significant level for $H_{0a}$ for 9
combinations of $(D_{2n},D_{3n})$, where the sizes corresponding
to $D_{2n}=0$ and \mbox{$D_{3n}=0$}. We observed that the size of the
ANOVA tests improved as the sample sizes and the observational
length $T$ increased, and the overall level of size were close to
the nominal $5\%$. This is quite reassuring considering the ANOVA
test is based on the asymptotic chi-square distribution. We also
observed that the power of the test increased as sample sizes and
$T$ were increased, and as the distance among the three
$\beta_{j0}$ was increased. For example, when $D_{2n}=0.0$ and
$D_{3n}=0.3$, the $L_2$ distance was $\sqrt{0.3^2+0.3^2}=0.424$,
which is larger than $\sqrt{0.1^2+0.2^2+0.3^2}=0.374$ for
$D_{2n}=0.2$ and $D_{3n}=0.3$. This explains why the ANOVA test
was more powerful for $D_{2n}=0.0$ and $D_{3n}=0.3$ than
$D_{2n}=0.2$ and $D_{3n}=0.3$. At the same time, we see similar
power performance between the two missing mechanisms.

%
\begin{table}[b]
\vspace*{-2pt}
\caption{Empirical size and power of the 5\% test for the existence of
interaction $H_{0c}\dvtx\gamma_{20}=0$}
\label{table2}
\begin{tabular*}{\tablewidth}{@{\extracolsep{\fill}}ccclcccccc@{}}
\hline
\multicolumn{3}{@{}c}{\textbf{Sample size}} & & &
\multicolumn{2}{c}{\textbf{Missingness}}
& & \multicolumn{2}{c@{}}{\textbf{Missingness}}\\[-4pt]
\multicolumn{3}{@{}c}{\hrulefill} & & & \multicolumn{2}{c}{\hrulefill}
& & \multicolumn{2}{c@{}}{\hrulefill}\\
$\bolds{n_1}$ & $\bolds{n_2}$ & $\bolds{n_3}$
& \multicolumn{1}{c}{$\bolds{\gamma_{20}}$} & $\bolds{T}$ & \textbf{I}
& \textbf{II} & $\bolds{T}$ & \textbf{I} & \textbf{II}\\
\hline
\phantom{0}60 & \phantom{0}65 & \phantom{0}55 & 0.0 (size) & 5 & 0.052 & 0.048 & 10 & 0.048 & 0.052 \\
& & & 0.2 & & 0.428 & 0.456 & & 0.568 & 0.636 \\
& & & 0.3 & & 0.722 & 0.788 & & 0.848 & 0.882 \\
& & & 0.4 & & 0.928 & 0.952 & & 0.948 & 0.968 \\
100 & 110 & 105 & 0.0 (size) & 5 & 0.054 & 0.046 & 10 & 0.056 & 0.042\\
& & & 0.2 & & 0.608 & 0.718 & & 0.694 & 0.812 \\
& & & 0.3 & & 0.940 & 0.938 & & 0.940 & 0.958 \\
& & & 0.4 & & 0.986 & 0.994 & & 0.952 & 0.966 \\
\hline
\end{tabular*}
\end{table}

To gain information on the empirical performance of the test on the existence
of interaction, we carried out a test for
$H_{0c}\dvtx\gamma_{20}=0$. In the simulation, we chose
$\gamma_{20}=0,0.2,0.3,0.4$, $\beta_{20}=2+\gamma_{20}$ and fixed
$\Delta_{2n}(t)=0$, respectively.\vadjust{\goodbreak} Table~\ref{table2} summarizes the sizes and
the powers of the test. Table~\ref{table3} reports the simulation
results of the ANOVA test on the interaction effect
$H_{0d}\dvtx\gamma_{10}=\gamma_{20}=\gamma_{30}$ with a similar
configurations as those
used as the ANOVA tests for the covarites effects reported in Table~\ref{table1}.
We observe satisfactory performance of these two tests in terms
of both the accurate of the size approximation and the empirical power.
In particular, the performance of the ANOVA tests were very much
similar to that conveyed in
Table~\ref{table1}.

%
\begin{table}
\caption{Empirical size and power of the 5\% ANOVA test for
$H_{0d}\dvtx\gamma_{10}=\gamma_{20}=\gamma_{30}$}
\label{table3}
\begin{tabular*}{\tablewidth}{@{\extracolsep{\fill}}cccclcccccc@{}}
\hline
\multicolumn{3}{@{}c}{\textbf{Sample size}} & & & & \multicolumn{2}{c}{\textbf{Missingness}} & &
\multicolumn{2}{c@{}}{\textbf{Missingness}}\\[-4pt]
\multicolumn{3}{@{}c}{\hrulefill} & & & & \multicolumn{2}{c}{\hrulefill} & &
\multicolumn{2}{c@{}}{\hrulefill}\\
$\bolds{n_1}$ & $\bolds{n_2}$ & $\bolds{n_3}$ & $\bolds{D_{2n}}$
& \multicolumn{1}{c}{$\bolds{D_{3n}}$} & $\bolds{T}$ & \textbf{I} & \textbf{II} & $\bolds{T}$
& \textbf{I} & \textbf{II}\\
\hline
\phantom{0}60 & \phantom{0}65 & \phantom{0}55 & 0.0 & 0.0 (size) & 5 & 0.058 & 0.058 & 10 & 0.068 &
0.036 \\
& & & 0.2 & 0.0 & & 0.134 & 0.188 & & 0.232 & 0.254 \\
& & & 0.3 & 0.0 & & 0.358 & 0.486 & & 0.510 & 0.622 \\
& & & 0.0 & 0.2 & & 0.136 & 0.166 & & 0.230 & 0.218 \\
& & & 0.0 & 0.3 & & 0.356 & 0.414 & & 0.466 & 0.474 \\
& & & 0.2 & 0.2 & & 0.170 & 0.208 & & 0.286 & 0.276 \\
& & & 0.2 & 0.3 & & 0.292 & 0.328 & & 0.462 & 0.428 \\
& & & 0.3 & 0.2 & & 0.266 & 0.356 & & 0.498 & 0.474 \\
& & & 0.3 & 0.3 & & 0.392 & 0.476 & & 0.578 & 0.588 \\
100 & 110 & 105 & 0.0 & 0.0 (size) & 5 & 0.068 & 0.040 & 10 & 0.054 &
0.046\\
& & & 0.2 & 0.0 & & 0.262 & 0.366 & & 0.354 & 0.432 \\
& & & 0.3 & 0.0 & & 0.654 & 0.744 & & 0.744 & 0.820 \\
& & & 0.0 & 0.2 & & 0.272 & 0.330 & & 0.340 & 0.334 \\
& & & 0.0 & 0.3 & & 0.590 & 0.676 & & 0.722 & 0.672 \\
& & & 0.2 & 0.2 & & 0.282 & 0.332 & & 0.412 & 0.410 \\
& & & 0.2 & 0.3 & & 0.528 & 0.582 & & 0.716 & 0.640 \\
& & & 0.3 & 0.2 & & 0.502 & 0.580 & & 0.680 & 0.728 \\
& & & 0.3 & 0.3 & & 0.672 & 0.674 & & 0.814 & 0.808 \\
\hline
\end{tabular*}
\end{table}

We then evaluate the power and size of the proposed ANOVA test
regarding the nonparametric components. To study the local power
of the test, we set $\Delta_{2n}(t)=U_n\sin(2\pi t)$ and
$\Delta_{3n}(t)=2\sin(2\pi t)-2\sin(2\pi(t+V_n))$, and fixed
$D_{2n}=0$ and $D_{3n}=0.2$. Here, $U_{n}$ and $V_n$ were designed
to adjust the amplitude and phase of the sine function. The same
kernel and bandwidths chosen by the cross-validation as outlined
earlier in the parametric ANOVA test were used in the test for
the nonparametric time effects. We calculated the test statistic
${\mathcal{T}}_n$ with $\varpi(t)$ being the kernel density estimate
based on all the time points in all treatments. We applied the
wild bootstrap proposed in Section~\ref{sec6} with $B=100$ to obtain
$\hat{t}_{0.05}$, the bootstrap estimator of the 5\% critical
value. The simulation results of the nonparametric ANOVA test for
the time effects are given in Table~\ref{table4}.

The sizes of the nonparametric ANOVA test were obtained when
$U_n=0$ and $V_n=0$, which were quite close to the nominal 5\%. We
observe that the power of the test increased when the distance
among $g_1(\cdot)$, $g_2(\cdot)$ and $g_3(\cdot)$ were becoming
larger, and when the sample size or repeated measurement $T$ were
increased. We noticed that the power was more sensitive to change
in $V_n$, the initial phase of the sine function, than $U_n$.

We then compared the proposed tests with
a test proposed by \citet{r34}.
Scheike and Zhang's test was comparing two treatments for the
nonparametric regression model
(\ref{eq:model_nonpara}) for longitudinal data without missing
values. Their test was based on a cumulative statistic
\[
T(z)=\int_a^z\bigl(\hat{g}_1(t)-\hat{g}_2(t)\bigr)\,dt,
\]
where $a,z$ are in a common time interval $[0, 1]$. They showed
that $\sqrt{n_1+n_2}T(z)$ converges to a Gaussian Martingale with
mean 0 and variance function
$\rho_1^{-1}h_1(z)+\rho_2^{-1}h_2(z)$, where
$h_j(z)=\int_a^z\sigma^2_{j}(y)f^{-1}_j(y)\,dy$. Hence, the test
statistic $T(1-a)/\sqrt{\widehat{\operatorname{Var}}\{T(1-a)\}}$ is used
for two group time-effect functions comparison.\vspace*{2pt}

%
\begin{table}
\caption{Empirical size and power of the 5\% ANOVA test for
$H_{0b}\dvtx g_{1}(\cdot)=g_{2}(\cdot)=g_{3}(\cdot)$ with
$\Delta_{2n}(t)=U_n\sin(2\pi t)$ and $\Delta_{3n}(t)=2\sin(2\pi
t)-2\sin(2\pi(t+V_n))$}\vspace*{-3pt}
\label{table4}
\begin{tabular*}{\tablewidth}{@{\extracolsep{\fill}}lccclcccccc@{}}
\hline
\multicolumn{3}{@{}c}{\textbf{Sample size}} & & & & \multicolumn{2}{c}{\textbf{Missingness}} & &
\multicolumn{2}{c@{}}{\textbf{Missingness}}\\[-4pt]
\multicolumn{3}{@{}c}{\hrulefill} & & & & \multicolumn{2}{c}{\hrulefill} & &
\multicolumn{2}{c@{}}{\hrulefill}\\
$\bolds{n_1}$ & $\bolds{n_2}$ & $\bolds{n_3}$ & $\bolds{U_n}$
& \multicolumn{1}{c}{$\bolds{V_n}$} & $\bolds{T}$ & \textbf{I} & \textbf{II} & $\bolds{T}$
& \textbf{I} & \textbf{II}\\
\hline
\phantom{0}60 & \phantom{0}65 & \phantom{0}55 & 0.00 & 0.00 (size) & 5 & 0.040 & 0.050 & 10 & 0.054 &
0.060 \\
& & & 0.30 & 0.00 & & 0.186 & 0.232 & & 0.282 & 0.256\\
& & & 0.50 & 0.00 & & 0.666 & 0.718 & & 0.828 & 0.840\\
& & & 0.00 & 0.05 & & 0.664 & 0.726 & & 0.848 & 0.842\\
& & & 0.00 & 0.10 & & 1.000 & 1.000 & & 1.000 & 1.000\\
100 & 110 & 105 & 0.00 & 0.00 (size) & 5 & 0.032 & 0.062 & 10 & 0.050 &
0.036\\
& & & 0.30 & 0.00 & & 0.434 & 0.518 & & 0.526 & 0.540\\
& & & 0.50 & 0.00 & & 0.938 & 0.980 & & 0.992 & 0.998\\
& & & 0.00 & 0.05 & & 0.916 & 0.974 & & 1.000 & 1.000\\
& & & 0.00 & 0.10 & & 1.000 & 1.000 & & 1.000 & 1.000\\
\hline
\end{tabular*}
\end{table}

To make the proposed test and the test of \citet{r34} comparable,
we conducted simulation in a set-up that mimics
the setting of model (\ref{eq5.4}) but with only the first two
treatments, no missing values and only the nonparametric part in
the regression by setting $\beta_j=\gamma_j=0$. Specifically, we
test for $H_0\dvtx g_1(\cdot) = g_2(\cdot)$ vs. $H_1\dvtx g_1(\cdot) =
g_2(\cdot) + \Delta_{2 n}(\cdot)$ for three cases of the
alternative shift function $\Delta_{2n}(\cdot)$ functions which
are spelt out in Table~\ref{table5} and set $a=0$ in the test of Scheike and
Zhang.
%
%
\begin{table}
\caption{The empirical sizes and powers of the proposed test (CZ)
and the test (SZ) proposed by Scheike and Zhang (\protect\citeyear{r34})
for $H_{0b}\dvtx g_{1}(\cdot)=g_{2}(\cdot)$ vs.
$H_{1b}\dvtx g_{1}(\cdot)=g_{2}(\cdot) + \Delta_{2n}(\cdot)$}
\label{table5}
\begin{tabular*}{\tablewidth}{@{\extracolsep{\fill}}lcclcccccc@{}}
\hline
\multicolumn{3}{@{}c}{\textbf{Sample size}} & & & \multicolumn{2}{c}{\textbf{Tests}}
& & \multicolumn{2}{c@{}}{\textbf{Tests}}\\[-4pt]
\multicolumn{3}{@{}c}{\hrulefill} & & & \multicolumn{2}{c}{\hrulefill}
& & \multicolumn{2}{c@{}}{\hrulefill}\\
$\bolds{n_1}$ & $\bolds{n_2}$ & $\bolds{n_3}$ &
\multicolumn{1}{c}{$\bolds{U_n}$}
& $\bolds{T}$ & \textbf{CZ} & \textbf{SZ} & $\bolds{T}$ & \textbf{CZ}
& \textbf{SZ}\\
\hline
\phantom{0}60 & \phantom{0}65 & \phantom{0}55 & \multicolumn{7}{c@{}}{Case I: $\Delta_{2n}(t)=U_n\sin
(2\pi t)$} \\[2pt]
& & & 0.00 (size) & 5 & 0.060 & 0.032 & 10 & 0.056 & 0.028 \\
& & & 0.30 & & 0.736 & 0.046 & & 0.844 & 0.028\\
& & & 0.50 & & 1.000 & 0.048 & & 1.000 & 0.026\\[2pt]
& & & \multicolumn{7}{c@{}}{Case II: $\Delta_{2n}(t)=2\sin(2\pi
t)-2\sin(2\pi(t+U_n))$} \\[2pt]
& & & 0.05 & & 1.000 & 0.026 & & 1.000 & 0.042\\
& & & 0.10 & & 1.000 & 0.024 & & 1.000 & 0.044\\[2pt]
& & & \multicolumn{7}{c@{}}{Case III: $\Delta_{2n}(t)=-U_n$} \\[2pt]
& & & 0.10 & & 0.196 & 0.162 & & 0.206 & 0.144\\
& & & 0.20 & & 0.562 & 0.514 & & 0.616 & 0.532
\\[2pt]
100& 110& 105 & \multicolumn{7}{c@{}}{Case I: $\Delta_{2n}(t)=U_n\sin
(2\pi t)$} \\[2pt]
& & & 0.00 (size) & 5 & 0.056 & 0.028 & 10 & 0.042 & 0.018 \\
& & & 0.30 & & 0.982 & 0.038 & & 0.994 & 0.040\\
& & & 0.50 & & 1.000 & 0.054 & & 1.000 & 0.028\\[2pt]
& & & \multicolumn{7}{c@{}}{Case II: $\Delta_{2n}(t)=2\sin(2\pi
t)-2\sin(2\pi(t+U_n))$} \\[2pt]
& & & 0.05 & & 1.000 & 0.022 & & 1.000 & 0.030\\
& & & 0.10 & & 1.000 & 0.026 & & 1.000 & 0.030\\[2pt]
& & & \multicolumn{7}{c@{}}{Case III: $\Delta_{2n}(t)=-U_n$}
\\[2pt]
& & & 0.10 & & 0.290 & 0.260 & & 0.294 & 0.218 \\
& & & 0.20 & & 0.780 & 0.774 & & 0.760 & 0.730 \\
\hline
\end{tabular*}
\end{table}
The simulation results are summarized in Table~\ref{table5}. We
found that in the first two cases (I and II) of the alternative shift
function $\Delta_{2 n}$, the test of
Scheike and Zhang had little power. It was only in the third case
(III), the test started to pick up some power although it was still not
as powerful as the proposed test.

\section{Analysis on HIV-CD4 data}\label{sec8}

In this section, we analyzed a longitudinal data set from AIDS
Clinical Trial Group 193A Study [\citet{r13}],
which was a randomized, double-blind study of HIV-AIDS patients
with advanced immune suppression. The study was carried out in
1993 with 1309 patients who were randomized to four treatments
with regard to HIV-1 reverse transcriptase inhibitors. Patients
were randomly assigned to one of four daily treatment regimes:
600 mg of zidovudine alternating monthly with 400 mg didanosine
(treatment I); 600 mg of zidovudine plus 2.25 mg of zalcitabine
(treatment II); 600 mg of zidovudine plus 400 mg of didanosine
(treatment III); or 600 mg of zidovudine plus 400~mg of didanosine
plus 400 mg of nevirapine (treatment VI). The four treatments had
325, 324, 330 and 330 patients, respectively.

The aim of our analysis was to compare the effects of age (Age),
baseline CD4 counts (PreCD4) and gender (Gender) on $Y =$ log(CD4
counts $+$1). The semiparametric model regression is, for $j=1,2,3$
and $4$,
%
%
\begin{equation}\label{model-realdata}
\quad
Y_{ji}(t)=\beta_{j1}
\operatorname{Age}_{ji}+\beta_{j2} \operatorname{PreCD4}_{ji} +\beta_{j3}
\operatorname{Gender}_{ji} +g_j(t)+\varepsilon_{ji}(t)
\end{equation}
with the intercepts absorbed in the nonparametric $g_j(\cdot)$
functions, and $\beta_j=(\beta_{j1},\beta_{j2},\beta_{j3})^\tau$ is the
regression coefficients to the covariates (Age, PreCD4, Gender).

To make $g_j(t)$ more interpretable, we centralized Age and
PreCD4 so that their sample means in each treatment were 0,
respectively. As a result, $g_j(t)$ can be interpreted as the
baseline evolution of $Y$ for a female (Gender${}={}$0) with the average
PreCD4 counts and the average age in treatment $j$. This kind of
normalization is used in \citet{r44} in their analyzes
for another CD4 data set. Our objectives were to detect any
difference in the treatments with respect to (i) the covariates;
and (ii) the nonparametric baseline functions.

Measurements of CD4 counts were scheduled at the start time 1 and
at a 8-week intervals during the follow-up. However, the data were
unbalanced due to variations from the planned measurement time and
missing values resulted from skipped visits and dropouts. The
number of CD4 measurements for patients during the first 40 weeks
of follow-up varied from 1 to 9, with a median of 4. There were
5036 complete measurements of CD4, and 2826 scheduled measurements
were missing. Hence, considering missing values is very important
in this analysis. Most of the missing values follow the monotone
pattern. Therefore, we model the missing mechanism under the
monotone assumption.

We considered three logistic regression models for the missing
propensities and used the AIC and BIC criteria to select the one
that was the mostly supported by data. The first model (M1) was a
logistic regression model for
$p_{j}(\accentset{\leftharpoonup}{X}_{jit,3},\accentset{\leftharpoonup
}{Y}_{jit,3};\theta_{j0})$
that effectively depends on $X_{jit}$ (the PreCD4) and
$(Y_{ji(t-1)}$, $Y_{ji(t-2)},Y_{ji(t-3)})$
if $t>3$. For $t<3$, it relies on all $Y_{jit}$
observed before $t$. In the second model (M2), we replace the
$X_{jit}$ in the first model with an intercept. In the third model
(M3), we added to the second logistic model with covariates
representing the square of $Y_{ji (t-1)}$ and the interactions
between $Y_{ji (t-1)}$ and $Y_{ji (t-2)}$. In the formulation of
the AIC and BIC criteria, we used the binary conditional likelihood
given in (\ref{BinLikehood}) with the respective penalties. The
difference of AIC and BIC values among these models for four
treatment groups is given in Table~\ref{table6}. Under the BIC criterion, M2
was the best model for all four treatments. For treatments II and
III, M3 had smaller AIC values than M2, but the differences were
very small. For treatments I and VI, M2 had smaller AIC than M3.
As the AIC tends to select more explanatory variables, we chose
M2 as the model
for the parametric missing propensity. 

%
\begin{table}
\def\arraystretch{0.9}
\caption{Difference in the AIC and BIC scores among three models
(\textup{M1})--(\textup{M3})}
\label{table6}
\begin{tabular*}{\tablewidth}{@{\extracolsep{\fill}}ld{2.2}d{3.2}d{2.2}d{2.2}d{2.2}d{2.2}d{2.2}d{3.2}@{}}
\hline
& \multicolumn{2}{c}{\textbf{Treatment I}} & \multicolumn{2}{c}{\textbf{Treatment
II}}
& \multicolumn{2}{c}{\textbf{Treatment III}} & \multicolumn{2}{c@{}}{\textbf{Treatment
VI}}\\[-4pt]
& \multicolumn{2}{c}{\hspace*{-1.5pt}\hrulefill} & \multicolumn{2}{c}{\hrulefill}
& \multicolumn{2}{c}{\hrulefill} & \multicolumn{2}{c@{}}{\hspace*{-1.5pt}\hrulefill}\\
\textbf{Models} & \multicolumn{1}{c}{\textbf{AIC}} & \multicolumn{1}{c}{\textbf{BIC}}
& \multicolumn{1}{c}{\textbf{AIC}} & \multicolumn{1}{c}{\textbf{BIC}}
& \multicolumn{1}{c}{\textbf{AIC}} & \multicolumn{1}{c}{\textbf{BIC}}
& \multicolumn{1}{c}{\textbf{AIC}} & \multicolumn{1}{c@{}}{\textbf{BIC}} \\
\hline
(M1)-(M2) & 3.85 & 3.85 & 14.90 & 14.90 & 17.91 & 17.91 & 10.35 & 10.35\\
(M2)-(M3) & -2.47& -11.47 & 0.93 & -8.12 & 0.30 & -8.75 & -3.15 &-12.27\\
\hline
\end{tabular*}
\end{table}

Model (\ref{model-realdata}) does not have interactions. It is
interesting to check if there is an interaction between gender and time.
Then the model becomes
%
%
\begin{eqnarray}
\label{model-realdata1}
Y_{ji}(t)&=&\beta_{j1}
\mathrm{Age}_{ji}+{\beta_{j2} \operatorname{PreCD4}_{ji}} +
{\beta_{j3}
\operatorname{Gender}_{ji}}\nonumber\\[-8pt]\\[-8pt]
&&{} +{\gamma_{j4}\operatorname
{Gender}_{ji}}\times t
+g_j(t)+\varepsilon_{ji}(t).\nonumber
\end{eqnarray}
We applied the proposed test in
Section~\ref{sec5} for $H_{0c}\dvtx\gamma_{j4}=0$ for $j=1,2,3$ and~$4$,
respectively. The
$p$-values were $0.9234,0.9885,0.9862$ and 0.5558, respectively, which means
that the interaction was not significant. Therefore, in the following
analyzes, we would not include the interaction term and continue
to use model (\ref{model-realdata}).

Table~\ref{table7} reports the parameter estimates $\hat{\beta}_j$ of
$\beta_{j}$ based on the estimating function $Z_{ji}(\beta_j)$
given in Section~\ref{sec3}. It contains the standard errors of the
estimates, which were obtained from the length of the EL
confidence intervals based on the marginal empirical likelihood
ratio for each $\beta_j$ as proposed in \citet{r2}. In
getting these estimates, we use the ``leave-one-subject''
cross-validation [\citet{r31}] to select the smoothing
bandwidths $\{h_j\}_{j=1}^4$ for the four treatments, which were
$12.90, 7.61, 8.27$ and $16.20$, respectively. We see that the
estimates of the coefficients for the Age and PreCD4 were similar
among all four treatments with comparable standard errors,
respectively. In particular, the estimates of the Age coefficients
endured large variations while the estimates of the PreCD4
coefficients were quite accurate. However, estimates of the Gender
coefficients had different signs among the treatments. We may also
notice that the confidence intervals from treatments I--IV for each
coefficient were overlap.

%
\begin{table}[b]
\def\arraystretch{0.9}
\caption{Parameter estimates and their standard errors}
\label{table7}
\begin{tabular*}{\tablewidth}{@{\extracolsep{\fill}}lcccc@{}}
\hline
& \textbf{Treatment I} & \textbf{Treatment II} & \textbf{Treatment III}
& \textbf{Treatment IV}\\[-4pt]
& \multicolumn{4}{c@{}}{\hrulefill}\\
\textbf{Coefficients} & $\bolds{\beta_1}$ & $\bolds{\beta_2}$ & $\bolds{\beta_3}$
& $\bolds{\beta_4}$ \\
\hline
Age & 0.0063\ (0.0039) & 0.0050\ (0.0040) & \phantom{$-$}0.0047\ (0.0058) & \phantom{$-$}0.0056\ (0.0046)
\\
PreCD4 & 0.7308\ (0.0462) & 0.7724\ (0.0378) & \phantom{$-$}0.7587\ (0.0523) & \phantom{$-$}0.8431\ (0.0425)
\\
Gender & 0.1009\ (0.0925) & 0.1045\ (0.0920) & $-$0.3300\ (0.1510) & $-$0.3055\ (0.1136)
\\
\hline
\end{tabular*}
\end{table}

We then formally tested $H_{0a}\dvtx\beta_1=\beta_2=\beta_3=\beta_4$.
The empirical likelihood ratio statistic $\ell_n$ was 8.1348,
which was smaller than $\chi_{9,0.95}^2=16.9190$, which produced a
$p$-value of 0.5206. So we do not have enough evidence to reject
$H_{0a}$ at a significant level 5 \%. The parameter estimates
reported in Table~\ref{table7} suggested similar covariate effects between
treatments I and II, and between treatments III and IV,
respectively; but different effects between the first two
treatments and the last two treatments. To verify this suggestion,
we carry out formal ANOVA test for pair-wise equality among the
$\beta_j$'s as well as for equality of any three\vadjust{\goodbreak} $\beta_j$'s. The
$p$-values of these ANOVA test are reported in Table~\ref{table8}. Indeed, the
difference between the first two treatments and between the last
two treatments were insignificant. However, the differences
between the first three (I, II and III) treatments and the last
treatment were also not significant.

%
\begin{table}
\caption{$p$-values of ANOVA tests for $\beta_j$'s}
\label{table8}
\begin{tabular*}{\tablewidth}{@{\extracolsep{\fill}}lccc@{}}
\hline
$\bolds{H_{0a}}$ & $\bolds{p}$\textbf{-value} & $\bolds{H_{0a}}$ & $\bolds{p}$\textbf{-value} \\
\hline
$\beta_1=\beta_2$& 0.9661 & $\beta_1=\beta_2=\beta_3$ & 0.7399\\
$\beta_1=\beta_3$& 0.4488 & $\beta_1=\beta_2=\beta_4$ & 0.4011\\
$\beta_1=\beta_4$& 0.1642 & $\beta_1=\beta_3=\beta_4$ & 0.3846\\
$\beta_2=\beta_3$& 0.4332 & $\beta_2=\beta_3=\beta_4$ & 0.4904\\
$\beta_2=\beta_4$& 0.2523 & $\beta_1=\beta_2=\beta_3=\beta_4$ &
0.5206 \\
$\beta_3=\beta_4$& 0.8450 & & \\
\hline
\end{tabular*}
\end{table}

%
\begin{figure}[b]

\includegraphics{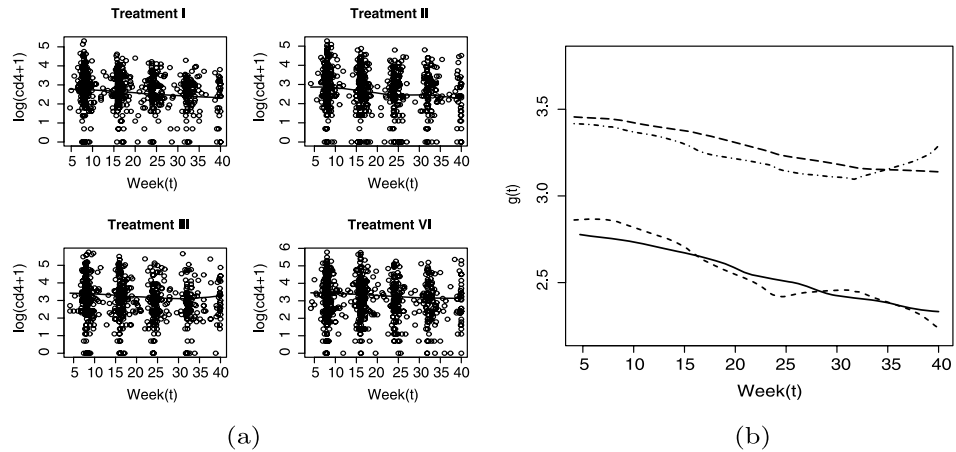}

\caption{\textup{(a)} The raw data excluding missing values plots
with the estimates of $g_j(t)$ ($j=1,2,3,4$). \textup{(b)} The estimates of
$g_j(t)$ in the same plot: treatment \textup{I} (solid line), treatment
\textup{II} (short dashed line), treatment \textup{III} (dashed and doted line)
and treatment \textup{IV} (long dashed line).}\label{figure1}
\end{figure}

We then tested for the nonparametric baseline time effects. The
kernel estimates $\hat{g}_j(t)$ are displayed in Figure~\ref{figure1}, which
shows that treatments I and II and treatments~III and IV had
similar baselines evolution overtime, respectively. However, a big
difference existed between the first two treatments and the last
two treatments. Treatment IV decreased more slowly than that of
the other three treatments, which seemed to be the most effective
in slowing down the decline of CD4. We also found that during the
first 16 weeks the CD4 counts decrease slowly and then the decline
became faster after 16 weeks for treatments I, II and III.

The $p$-value for testing $H_{0b}\dvtx
g_1(\cdot)=g_2(\cdot)=g_3(\cdot)=g_4(\cdot)$ is shown in Table~\ref{table9}.
The entries were based on 500 bootstrapped resamples according to
the procedure introduced in Section~\ref{sec6}.
The statistics ${\mathcal{T}}_n$ for testing $H_{0b}\dvtx
g_1(\cdot)=g_2(\cdot)=g_3(\cdot)=g_4(\cdot)$ was 3965.00, where we
take $\varpi(t)=1$ over the range of $t$. The $p$-value of the test
was 0.004. Thus, there existed significant difference in the
baseline time effects $g_j(\cdot)$'s among treatments I--IV. At the
same time, we also calculate the test statistics ${\mathcal{T}}_n$ for
testing $g_1(\cdot)=g_2(\cdot)$ and $g_3(\cdot)=g_4(\cdot)$. The
statistics values were 19.26 and 26.22, with $p$-values 0.894 and
0.860, respectively. These $p$-values are much bigger than 0.05. We
conclude that treatment I and II has similar baseline time
effects, but they are significantly distinct from the baseline
time effects of treatment III and IV, respectively. $p$-values of
testing other combinations on equalities of $g_1(\cdot),
g_2(\cdot), g_3(\cdot)$ and $g_4(\cdot)$ are also reported in
Table~\ref{table9}.

%
\begin{table}
\caption{$p$-values of ANOVA tests on $g_j(\cdot)$'s}
\label{table9}
\begin{tabular*}{\tablewidth}{@{\extracolsep{\fill}}lccc@{}}
\hline
$\bolds{H_{0b}}$ & $\bolds{p}$\textbf{-value} & $\bolds{H_{0b}}$
& $\bolds{p}$\textbf{-value} \\
\hline
$g_1(\cdot)=g_2(\cdot)$& 0.894 & $g_1(\cdot)=g_2(\cdot)=g_3(\cdot
)$ & 0.046\\
$g_1(\cdot)=g_3(\cdot)$& 0.018 & $g_1(\cdot)=g_2(\cdot)=g_4(\cdot
)$ & 0.010\\
$g_1(\cdot)=g_4(\cdot)$& 0.004 & $g_1(\cdot)=g_3(\cdot)=g_4(\cdot
)$ & 0.000\\
$g_2(\cdot)=g_3(\cdot)$& 0.020 & $g_2(\cdot)=g_3(\cdot)=g_4(\cdot
)$ & 0.014\\
$g_2(\cdot)=g_4(\cdot)$& 0.006 & $g_1(\cdot)=g_2(\cdot)=g_3(\cdot
)=g_4(\cdot)$ &
0.004 \\
$g_3(\cdot)=g_4(\cdot)$& 0.860 & & \\
\hline
\end{tabular*}
\end{table}

This data set has been analyzed by \citet{r9} using a random
effects model that applied the Restricted
Maximum Likelihood (REML) method. They conducted a two sample
comparison test via parameters in the model for the difference
between the dual therapy (treatment I--III) versus triple therapy
(treatment VI) without considering the missing values. More
specifically, they denoted $\mathrm{Group}=1$ if subject in the
triple therapy treatment and $\mathrm{Group}=0$ if subject in the
dual therapy treatment, and the linear mixed effect was
\begin{eqnarray*}
E(Y|b)&=&\beta_1+\beta_2t+\beta_3(t-16)_{+}+{\beta_4
\operatorname{Group}}\times t\\
&&{} +{\beta_5\operatorname{Group}}\times(t-16)_{+}+b_1+b_2t+b_3(t-16)_{+},
\end{eqnarray*}
where $b=(b_1,b_2,b_3)$ are random effects. They tested $H_0\dvtx
\beta_4=\beta_5=0$. This is equivalent to test the null hypothesis
of no treatment group difference in the changes in $\log$ CD4
counts between therapy and dual treatments. Both Wald test and
likelihood ratio test rejected the null hypothesis, indicating the
difference between dual and triple therapy in the change of $\log$
CD4 counts. Their results are consistent with the result we
illustrated in Table~\ref{table9}.\looseness=1

\begin{appendix}\label{app}
\section*{Appendix: Technical assumptions}

We provides the conditions used for Theorems~\ref{th1}--\ref{th4} and
some remark in this section. The proofs for Theorems
\ref{th1},~\ref{th2},~\ref{th3} and~\ref{th4} are contained in the supplement article
[\citet{r4}]. The proof for Theorem~\ref{th5} is largely similar
to that of Theorem~\ref{th1} and is omitted.

The following assumptions are made in the paper:

\begin{enumerate}[A1.]
\item[A1.] Let $S(\theta_j)$ be the score function of the partial
likelihood ${\mathcal{L}}_{B_j}(\theta_j)$ for a q-dimensional
parameter $\theta_j$ defined in (\ref{BinLikehood}), and
$\theta_{j0}$ is in the interior of compact $\Theta_j$. We assume
$E\{S(\theta_j)\}\neq0$ if $\theta_j\neq\theta_{j0}$,
$\operatorname{Var}(S(\theta_{j0}))$ is finite and positive definite, and
$E (\frac{\partial S(\theta_{j0})}{\partial
\theta_{j0}} )$ exists and is invertible. The missing
propensity $\pi_{jim}(\theta_{j0})>b_0>0$ for all $j,i,m$.\vspace*{1pt}

\item[A2.]
\begin{enumerate}[(iii)]
\item[(i)] The kernel function $K$ is a symmetric probability
density which is differentiable of Lipschitz order 1 on its
support $[-1,1]$. The bandwidths satisfy $n_jh_j^2/\log^2
n_j\rightarrow\infty$, $n_j^{1/2}h_j^4\rightarrow0$ and
$h_j\rightarrow0$ as $n_j\rightarrow\infty$.\vspace*{1pt}

\item[(ii)] For each treatment $j$ $(j=1,\ldots,k)$, the design points
$\{t_{jim}\}$ are thought to be independent and identically
distributed from a super-population with density $f_j(t)$. There
exist constants $b_l$ and $b_u$ such that $0<b_l \leq\sup_{t\in
S}f_j(t)\leq b_u<\infty$.

\item[(iii)] For each $h_j$ and $T_j$, $j=1,\ldots,k$, there exist
finite positive constants $\alpha_j$, $b_j$ and $T$ such that
$\alpha_j T_j=T$ and $b_j h_j = h$ for some $h$ as $h \to0$. Let
$n=\sum_{i=1}^k n_j, n_j/n\rightarrow\rho_j$ for some nonzero
$\rho_j$ as $n\rightarrow\infty$ such that $\sum_{i=1}^k\rho_j=1$.
\end{enumerate}

\item[A3.] The residuals $\{\varepsilon_{ji}\}$ and $\{u_{ji}\}$
are independent of each other and each of $\{\varepsilon_{ji}\}$
and $\{u_{ji}\}$ are mutually independent among different $j$ or
$i$,\break respectively; $\max_{1\leq i\leq
n_j}\|u_{jim}\|=o_p\{n_j^{({2+r})/({2(4+r)})}(\log n_j)^{-1}\}$,\break
$\max_{1\leq i\leq n_j}E|\varepsilon_{jim}|^{4+r}<\infty$, for
some $r>0$; and assume that
\[
\lim_{n_j\rightarrow\infty}(n_jT_j)^{-1}\sum_{i=1}^{n_j}\sum_{m=1}^{T_j}
E\{\widetilde{\mathbb{X}}_{jim}\widetilde{\mathbb{X}}_{jim}^\tau\}
=\Sigma_x>0,
\]
where
$\widetilde{\mathbb{X}}_{jim}=\mathbb{X}_{jim}-E(\mathbb{X}_{jim}|t_{jim})$.

\item[A4.] The functions $g_{j0}(t)$ and $h_j(t)$ are,
respectively, one-dimensional and $p$-dimensional smooth functions
with continuously second derivatives on $S=[0,1]$.
\end{enumerate}
\begin{remark*}
Condition A1 are the regular conditions for the
consistency of the binary MLE for the parameters in the missing
propensity. Condition A2(i) are the usual conditions for the
kernel and bandwidths in nonparametric curve estimation. Note that
the optimal rate for the bandwidth $h_j=O(n_j^{-1/5})$ satisfies
A2(i).\vadjust{\goodbreak} The requirement of design points $\{t_{jim}\}$ in A2(ii)
is a common assumption similar to the ones in \citet{r24}.
Condition A2(iii) is a mild assumption on the relationship between
bandwidths and sample sizes among different samples. In A3, we do
not require the residuals $\{\varepsilon_{ji}\}$ and $\{u_{ji}\}$
being, respectively, identically distributed for each fixed $j$.
This allows extra heterogeneity among individuals for a treatment.
The positive definite of $\Sigma_x$ in condition A3 is used to
identify the ``parameters'' $(\beta_{j0}, \gamma_{j 0}, g_{j 0})$
uniquely, which is a generalization of the identification
condition used in \citet{r11} to longitudinal
data. This condition can be checked empirically by constructing
consistent estimate of~$\Sigma_x$.
\end{remark*}
\end{appendix}

\section*{Acknowledgments}
The authors thank the referees,  Associate Editors and
Editors for valuable comments which lead to improvement of the
presentation of the paper.

\begin{supplement}
\stitle{Supplement to ``ANOVA for
Longitudinal Data with Missing Values''}
\slink[doi]{10.1214/10-AOS824SUPP}
\sdatatype{.pdf}
\sfilename{supplement.pdf}
\sdescription{This supplement material provides technical proofs
to the asymptotic distributions of the empirical likelihood ANOVA
test statistics for comparing the treatment effects with respect
to covariates given in Theorems~\ref{th1} and~\ref{th2}, the asymptotic
normality
of the empirical likelihood ratio based ANOVA test statistic for
comparing the nonparametric time effect functions given in Theorem
\ref{th3} and justifies the usage of the proposed bootstrap procedure.}
\end{supplement}


%
\printaddresses

\end{document}